\documentclass[final]{siamltex}

\usepackage{amsmath,amssymb}
\usepackage{showkeys}
\usepackage{psfrag}
\usepackage{pdfsync}
\usepackage{graphicx}
\usepackage{color}
\usepackage{hyperref}

\RequirePackage{snapshot}
\usepackage{ifpdf}
\ifpdf
\usepackage{pdfsync}
\else
\usepackage[active]{srcltx}
\fi

\makeatletter
\renewcommand\theequation{\thesection.\arabic{equation}}
\@addtoreset{equation}{section}
\makeatother
  \makeatletter
  \def\proof{\@ifnextchar[\opargproof{\opargproof[\bf Proof \hfil\\ ]}}
  \def\opargproof[#1]{\par\noindent {\bf #1 }}
  
  \makeatother
\newtheorem{defi}{Definition}
\newtheorem{coro}{Corollary}
\newtheorem{question}{Question}

\def\Bt{\begin{thm}}
\def\Et{\end{thm}}
\def\Br{\begin{rem}}
\def\Er{\end{rem}}
\def\Bd{\begin{defi}}
\def\Ed{\end{defi}}
\def\Bc{\begin{coro}}
\def\Ec{\end{coro}}
\def\Bp{\begin{prop}}
\def\Ep{\end{prop}}
\def\Bl{\begin{lem}}
\def\El{\end{lem}}

\def\Be{\begin{equation*}}
\def\Ee{\end{equation*}}
\def\Ben{\begin{equation}}
\def\Een{\end{equation}}

\def\Ba{\begin{eqnarray*}}
\def\Ea{\end{eqnarray*}}
\def\Ban{\begin{eqnarray}}
\def\Ean{\end{eqnarray}}

\def\R{\mathbf{R}}
\def\ds{\displaystyle}

\def\p#1{p_{#1}}

\def\rij{{r_{i,j}}}

\def\lra#1{\left\langle #1 \right\rangle}

\def\H{\mathcal{H}}
\def\HH{\mathcal{H}_N}
\def\AA{\mathcal{S}_N}

\def\TT{\mathcal{T}}

\def\bfc{{\bf{c}}}
\def\bfe{{\bf{e}}}
\def\bff{{\bf{f}}}
\def\bfg{{\bf{g}}}
\def\bfn{{\bf{n}}}
\def\bfq{{\bf{q}}}
\def\bfr{{\bf{r}}}
\def\bft{{\bf{t}}}
\def\bfu{{\bf{u}}}
\def\bfv{{\bf{v}}}
\def\bfw{{\bf{w}}}
\def\bfx{{\bf{x}}}
\def\bfy{{\bf{y}}}
\def\bfz{{\bf{z}}}
\def\bfC{{\bf{C}}}
\def\bfF{{\bf{F}}}
\def\bfG{{\bf{G}}}
\def\bfK{{\bf{K}}}

\def\bfS{{\bf{S}}}
\def\bfT{{\bf{T}}}
\def\bfV{{\bf{V}}}
\def\bfW{{\bf{W}}}
\def\bfX{{\bf{X}}}
\def\bfOm{{\bf{\Omega}}}
\def\bfomega{{\boldsymbol \omega}}
\def\bfphi{{\boldsymbol \phi}}
\def\bfsig{{\boldsymbol \sigma}}
\def\bfId{{\bf{Id}}}

\DeclareMathOperator{\SO}{SO(3)}

\def\nab{\nabla}
\def\nabdot{\nabla \cdot}

\def\cross{\times}

\def\pO{{\partial \Omega}}
\def\pB{\partial B}

\def\VV{\bfV}
\def\VW{\bfW}
\def\FF{\bfF}
\def\Rot{\mathcal R}
\def\MetG{\mathcal G}
\def\Res{\mathbf{Res}}
\def\Re{\R}

\newcommand{\Luca}[1]{\textcolor{black}{#1}}
\newcommand{\Francois}[1]{\textcolor{black}{#1}}

\title{Optimally swimming Stokesian robots}

\author{Fran\c cois Alouges\thanks{CMAP UMR 7641, \'Ecole Polytechnique CNRS, Route de Saclay, 91128 Palaiseau Cedex - France ({\tt{francois.alouges@polytechnique.edu}})} \and 
Antonio DeSimone
\thanks{SISSA, International School of Advanced Studies, Via Bonomea 265, 34136 Trieste - Italy
 ({\tt{desimone@sissa.it}})}
 \and Luca Heltai\thanks{SISSA, International School of Advanced Studies, Via Bonomea 265, 34136 Trieste - Italy
 ({\tt{luca.heltai@sissa.it}})} \and Aline Lefebvre-Lepot\thanks{CMAP UMR 7641, \'Ecole Polytechnique CNRS, Route de Saclay, 91128 Palaiseau Cedex - France ({\tt{aline.lefebvre@polytechnique.edu}})} \and Beno\^\i t Merlet\thanks{CMAP UMR 7641, \'Ecole Polytechnique CNRS, Route de Saclay, 91128 Palaiseau Cedex - France ({\tt{benoit.merlet@polytechnique.edu}})}}
\begin{document}

\maketitle

\begin{abstract}
  We study self-propelled stokesian robots composed of assemblies of
  balls, in dimensions 2 and 3, and prove that they are able to
  control their position and orientation.  This is a result of
  \emph{controllability}, and its proof relies on applying Chow's
  theorem in an analytic framework, similar to what has been done in
  \cite{AlougesDeSimoneLefebvre-2007-a} for an axisymmetric system
  swimming along the axis of symmetry.  We generalize
  the analyticity result given in
  \cite{AlougesDeSimoneLefebvre-2007-a} to the situation
  where the swimmers can move either in a plane or in
  three-dimensional space, hence experiencing also rotations. We then
  focus our attention on energetically optimal strokes, which we are
  able to compute numerically. Some examples of computed optimal
  strokes are discussed in detail.
\end{abstract} 



\section{Introduction}

Self-propulsion at low Reynolds number is a problem of considerable
biological and biomedical relevance which has also great appeal from
the point of view of fundamental science.  Starting from the
pioneering work by Taylor \cite{Taylor-1951-a} and Lighthill
\cite{Lighthill-1952-a}, it has received a lot of attention in recent
years (see the encyclopedia article
\cite{AlougesDeSimoneLefebvre-2009-a} for an elementary introduction
and the review paper \cite{LaugaPowers-2009-a} for a comprehensive
list of references).

Both relevance for applications and theoretical interest stem from the
fact that one considers swimmers of small size. The Reynolds number
$Re=LV/\nu$ gives an estimate for the relative importance of intertial
to viscous forces for an object of size $L$ moving at speed $V$
through a newtonian fluid with kinematic viscosity $\nu$. Since in
applications $V$ rarely exceeds a few body lengths per second, if one
considers swimming in a given medium, say, water, then $Re$ is
entirely controlled by $L$. At small $L$, inertial forces are
negligible and, in order to move, micro-swimmers can only exploit the
viscous resistance of the surrounding fluid. The subtle consequences
of this fact (which are rather paradoxical when compared to the
intuition we can gain from our own swimming experience) are discussed
in \cite{Purcell-1977-a}.  For example, the motion of microswimmers is
geometric: the trajectory of a low $Re$ swimmer is entirely determined
by the sequence of shapes that the swimmer assumes. Doubling the rate
of shape changes simply doubles the speed at which the same trajectory
is traversed.  As observed in \cite{ShapereWilczek-1989-a}, this
suggests that there must be a natural, attractive mathematical
framework for this problem (which the authors, indeed, unveil).  On
the other hand, bacteria and unicellular organisms \emph{are} of
micron size, while artificial robots to be used non-invasively inside
human bodies for medical purposes \emph{must} be small.  Discovering
the secrets of biological micro-swimmers and controlling engineered
micro-robots requires a quantitative understanding of low $Re$
self-propulsion.

The basic problem of swimming is easy to state: given a (periodic)
time history of shapes of a swimmer (a sequence of strokes), determine
the corresponding time history of positions and orientations in
space. A natural, related question is the following: starting from a
given position and orientation, can the swimmer achieve any prescribed
position and orientation by performing a suitable sequence of strokes?
This is a question of \emph{controllability}.  The peculiarity of low
$Re$ swimming is that, since inertia is negligible, reciprocal shape
changes lead to no net motion, so the question of controllability may
become non trivial for swimmers that have only a few degrees of
freedom at their disposal to vary their shape. The well
\Francois{established} scallop theorem \cite{Purcell-1977-a} is
precisely a result of non-controllability.

Once controllability is known, i.e., it is shown that it is possible
to go from A to B, one can ask the question of how to go from A to B
at minimal energetic cost. This is a question of \emph{optimal
  control}.

In spite of the clear connections between low $Re$ self-propulsion and
control theory, this viewpoint has started to emerge only recently,
and mostly in the mathematical literature. Examples are
\cite{KoillerEhlersMontgomery-1996-a}, and the more recent
contributions \cite{Bressan-2008-a}, \cite{ChambrionMunnier2011}, \cite{Khapalov-2007-a} and
\cite{San-MartinTakahashiTucsnak-2007-a}.
The papers \cite{AlougesDeSimoneHeltai-2011-a},
\cite{AlougesDeSimoneLefebvre-2007-a},
\cite{AlougesDeSimoneLefebvre-2009-a} study in detail both
controllability and optimal control for axisymmetric swimmers whose
varying shapes are described by few (in particular, two) scalar
parameters.

In this paper we analyze the problem of low Reynolds number swimming
from the point of view of geometric control theory, and focus on a
special class of model systems: those obtained as assemblies of a
small number of balls.  Model systems of this type have played an important role in clarifying the subtleties of low $Re$ self-propulsion \cite{NajafiGolestanian-2004-a}, \cite{AvronGatKenneth-2004-a}, \cite{DreyfusBaudryStone-2005-a}. These systems offer an interesting balance
between complexity of the analysis and richness of observable
behavior. 

While constructing such artificial swimmers may indeed be possible,
  their practical uses are at present unclear. However, the analysis of their swimming patterns may prove a very useful tool both in the design of robotic microswimmers \cite{Dreyfusetal}, and in understanding the motion of biological swimmers. As an example, the social motility patterns exhibited by Myxobacteria (see, e.g.,~\cite{Harshey-2003-a}) are strikingly similar to the ones described in this paper. As part of their life cycle, individual cells become linked together and move collectively on surfaces. The links between adjacent cells are guaranteed by pili which individual cells can project and retract. The fact that collective motion rests upon control of the relative positions of individuals seems fully established. Our study may shed light on the many details that are at present unclear.
  
Proving controllability and providing optimal control strategies
  for realistic swimmers, either biological or artificial, is
  rather difficult. The use of assembly of balls as a model swimmer relies on
  the idea that one can use collections of spheres to model the
  hydrodynamic characteristics of real swimmers, replacing their
  actual swimming apparatus (which can be indeed very complex to
  describe, simulate and approximate) with a finite collection of spheres
  having comparable hydrodynamic resistance. Very recent progress towards the analysis of low $Re$ more general swimmers propelling through shape changes can be found in \cite{DalMasoDeSimoneMorandotti}, \cite{LoheacTucsnak}.

Our approach is similar in spirit to the one in
\cite{AlougesDeSimoneHeltai-2011-a},
\cite{AlougesDeSimoneLefebvre-2007-a},
\cite{AlougesDeSimoneLefebvre-2009-a}, but we extend it to
non-axisymmetric systems such as three spheres moving in a plane, and
systems of four spheres moving in three dimensional space. The motion of these systems
is described by both positional and orientational variables, leading
to a much richer geometric structure of the state space. The study of such systems requires substantial extensions of the mathematical
    and numerical methods introduced in
    \cite{AlougesDeSimoneLefebvre-2007-a}.

For all the model swimmers described above, controllability is proved
by using two main ingredients.  The first is Chow's theorem, leading
to local controllability in a neighborhood of a point $X$ in state
space. We verify the full rank~\eqref{fullrank} hypothesis of this
theorem by showing that the vector fields of the coefficients of the
governing ODEs (a linear control system without drift) and their first
order Lie brackets span the whole tangent space to the state space at
$X$.  The second is the analyticity of the coefficients, and the fact
that our shape space is always connected. This allows us to pass from local to global controllability generalizing a result that has been proved in \cite{AlougesDeSimoneLefebvre-2007-a} for the special axisymmetric case.  In the same
paper~\cite{AlougesDeSimoneLefebvre-2007-a}, Chow's theorem was also
used to prove local controllability. In this simpler axisymmetric
case, however, position is described only by one scalar parameter (the
position $p$ of one distinguished point along the axis of
symmetry). The non-degeneracy condition reduces then to the
non-vanishing of a single scalar quantity, namely, the curl of the
vector field $\VV$, governing the rate of change of position as a
consequence of shape changes at rates $\xi_i$ according to $\dot
p=\VV_1(\xi)\dot\xi_1 + \VV_2(\xi)\dot\xi_2$. In the richer context of
this paper, which involves also rotational degrees of freedom of the
swimmers, proving controllability requires an explicit computation of
all the first order Lie Brackets.  In fact, all the systems we analyze
here satisfy the condition
$$M+M(M-1)/2 = d$$
where $M$ is the number of controls (rate of change of shape
variables) and $d$ is the dimension of the state space (position,
orientation, and shape). This is the necessary condition that first
order Lie brackets alone suffice to show that the Lie algebra of the
coefficients has full rank, so that controllability follows.

For controllable systems, it makes sense to ask how
to achieve the desired target (position and orientation) at minimal
energy cost. We present a method to address this optimal control
question numerically, and we then examine in detail several
 optimal strokes for a concrete model swimmer (three balls swimming in a
plane with a prescribed lateral displacement). Depending on whether
the final orientation is also prescribed, and whether the initial
shape is prescribed as well or rather one treats it as a parameter to
be optimized, we obtain dramatically different answers. Their variety
illustrates the surprisingly richness of behavior of low $Re$
swimmers.

The rest of the paper is organized as follows.  \Francois{In Section
  \ref{sec1}, we describe the various model swimmers to which our analytical and numerical tools are later applied. The first is the
  Najafi-Golestanian's swimmer \cite{NajafiGolestanian-2004-a},
  already treated in \cite{AlougesDeSimoneLefebvre-2007-a}, while the
  two others are non-trivial generalizations. Section \ref{sec2}
  presents some results on Stokes flows and in Section \ref{sec3} we
  show that swimming is indeed an affine control problem
  without drift. In Section \ref{sec4}, we prove the effective
  swimming capability of our model swimmers. In Sections \ref{sec5} and
  \ref{sec6} we state the optimal control problem and a numerical
  strategy for its solution. Examples of optimal strokes for three
  balls swimming in a plane are discussed in detail in Section
  \ref{sec:numerical-results}.}

\section{The swimmers}
\label{sec1}
We will focus our attention on some special swimmers. Namely, 
we assume that the swimmer is composed of $N$ non-intersecting balls $(B_i)_{1\leq i \leq N}$ $\subset \mathbb{R}^3$
centered at $(\bfx_i)_{1\leq i \leq N}$, and we restrict ourselves 
to configurations which can be described by two sets of variables:
\begin{itemize}
\item the {\em{shape variables}}, denoted by $\xi \in \mathcal{S}$, where $\mathcal{S}$ is an open connected subset of $\R^M$, from which 
relative distances $(\bfx_{ij})_{1\leq i,j \leq N}$ between the balls $(B_i)_{1\leq i \leq N}$ are obtained\Francois{. In the examples treated in this paper, the balls are assumed to move only along fixed directions which make fixed angles one to another. This reflects a situation where the balls are linked together by thin jacks that are able to elongate. The viscous resistance associated with these jacks is, however, neglected, and the fluid is thus assumed to fill the whole set $\mathbb{R}^3\setminus \bigcup_{i=1}^N B_i$.}
\item the {\em{position variables}}, denoted by $p \in \mathcal{P}$, which
describe the global position and orientation in space of the swimmer. In our examples, the set $\mathcal{P}$ is typically a manifold of dimension less than or equal to six.
The six-dimensional case is given by  $p\in\mathcal{P}=\R^3\times \SO$ and $p$ consists of a translation and a rotation in the three-dimensional space.
\end{itemize}

We also assume that the orientation of the balls $(B_i)_{1\leq i \leq N}$ and the distances $(\bfx_{ij})_{1\leq i,j \leq N}$ depend analytically on $(\xi,p)$, therefore the
state of the system is analytically and uniquely determined by the variables $(\xi,p)$,
so that there exist $N$ (analytic) functions
\Be
\bfX_i : \mathcal{S} \times \mathcal{P} \times \pB \longrightarrow \R^3
\Ee
which give the position of the current point of the $i-$th sphere of the swimmer
in the state $(\xi,p)$ in $\mathcal{S}\times \mathcal{P}$. The non-slip boundary condition
on $\pB_i$ imposes that the velocity of the fluid is given by
\Ben
\bfv_i(\xi,p,\bfr)=\frac{d}{dt}\bfX_i(\xi,p,\bfr) = (\dot{\xi}\cdot \nabla_\xi) \bfX_i(\xi,p,\bfr)+(\dot{p}\cdot\nabla_p) \bfX_i(\xi,p,\bfr).
\label{eqvi}
\Een

\subsection{The three sphere swimmer of Najafi and Golestanian (3S)}
\label{swimmer1}
\begin{figure}
  \centering
  \includegraphics[width=.5\textwidth]{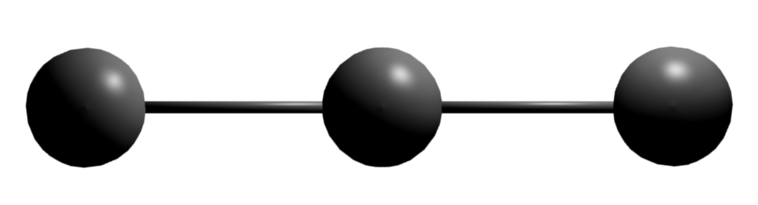}
  \caption{Three Sphere Swimmer (3S)}
  \label{NGswimmer}
\end{figure}

This swimmer, initially proposed in \cite{NajafiGolestanian-2004-a}  has been studied thoroughly 
in \cite{AlougesDeSimoneLefebvre-2007-a} and \cite{AlougesDeSimoneLefebvre-2008-a}. It is composed of 3 spheres of radius $a>0$ aligned along the $x-$axis,
as depicted in Fig. \ref{NGswimmer}.

We call $\xi=(\xi_1,\xi_2)$ the length of the arms, and $p$ the position of
the central sphere which leads us\Francois{, in order to avoid overlap of the spheres,} to consider the shape set $\mathcal{S}=(2a,+\infty)^2$ and the position set $\mathcal{P}=\R$. \Francois{Indeed, due to axial symmetry, this swimmer may only move along the $x-$axis.} Using $\bfe_1=(1,0,0)^T$, we can write
\Be 
\bfx_1(\xi,p)=(p-\xi_1)\bfe_1\,,
\bfx_2(\xi,p)=p\,\bfe_1\,,
\bfx_3(\xi,p)=(p+\xi_2)\bfe_1\,,
\Ee
and
\Be
\bfX_i(\xi,p,\bfr)=\bfx_i(\xi,p)+\bfr\,,\ \ \ \forall i\in\{1,2,3\}\,,\forall \bfr\in \partial B\,.
\Ee 

\subsection{The three sphere swimmer moving in a plane (3SP)}
\label{swimmer2}
\begin{figure}[!htb]
  \begin{center}

  \resizebox{.5\textwidth}{!}{
    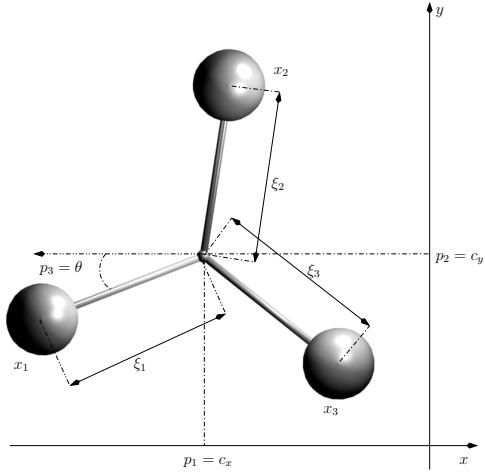
  }

    \caption{\label{3sub} The three-sphere planar swimmer (3SP)}
  \end{center}
\end{figure}

\Luca{A variant of the Najafi-Golestanian's swimmer, called
  ``\emph{Purcell's rotator}'', was presented
  in~\cite{DreyfusBaudryStone-2005-a}, where the axis along which the
  three spheres are constrained to move, is bent to form a circle with
  fixed radius. The resulting swimmer is one where the three spheres
  always keep a fixed distance from the center of the swimmer, and may
  only vary their location on the circle, remaining on the same plane.}

\Luca{As observed in~\cite{DreyfusBaudryStone-2005-a}, this variation
  allows one to control rotations around the axis perpendicular to the
  circle and passing through the center of the swimmer (hence the name
  ``\emph{rotator'}''), and introduces as well a small drift in the
  system, which displaces horizontally on the plane containing the
  circle.}

\Luca{That swimmer does not change the number of controls (two) with
  respect to the Najafi Golestanian's swimmer, however a periodic
  shape change induces both a rotation \emph{and} a translation of the
  swimmer (that is, the dimension $d$ of the state space is equal to
  $5 \neq M+M(M-1)/2 = 3$). Seen from the perspective of control
  theory, this implies that not only first order Lie brackets are used
  to change the positional variables but also higher order ones, with
  smaller and smaller effectiveness as the order of the used Lie
  brackets increases (see Section~\ref{sec4}).}

\Luca{A complementary version of Purcell's rotator, with one
  additional control variable,} has been proposed and studied in
\cite{Lefebvre-LepotMerlet-2009-a}. It is composed of three balls
$B_1,B_2,B_3$ of equal radii $a>0$.  The three balls can move along
three horizontal axes that mutually meet at $\bfc$ with angle $\ds
\frac{2\pi}{3}$ (see Fig.~\ref{3sub}). The balls do not rotate around
their axes so that the shape of the swimmer is characterized by the
three lengths $\xi_1,\xi_2,\xi_3$ of its arms, measured from the
origin to the center of each ball. However, the swimmer may freely
rotate around $\bfc$ in the horizontal plane.

Consider a reference equilateral triangle $(S_1,S_2,S_3)$ with
center $O \in \R^3$ in the horizontal plane $(O,x,y)$ such that
$\mbox{dist}(O,S_i)=1$ and define $\bft_i = \vec{O S}_i$.  
Position and orientation in the horizontal plane are
described by the coordinates of the center $\bfc\in \R^3$ (but $\bfc$
stays confined to the horizontal plane) and the horizontal angle $\theta$ that
one arm, say arm number 1, makes with a fixed direction, say
$(O,x)$, in such a way that $d=3$. Therefore, we place the center
of the ball $B_i$ at $\bfx_i=\bfc+\xi_i \Rot_\theta \bft_i$ with
$\xi_i>0$ for $i=1,2,3$, where $\Rot_\theta$ stands for the horizontal
rotation of angle $\theta$ given for instance by the matrix:
$$
\Rot_\theta=\left(
\begin{array}{ccc}
\cos(\theta) & -\sin(\theta) & 0\\
\sin(\theta) & \cos(\theta) & 0\\
0 & 0 & 1
\end{array}
\right)\,.
$$

The swimmer is then fully described by the parameters
$X=(\xi,\bfc,\theta)\in \mathcal{S}\times \mathcal{P}$, where $\mathcal{S} \ :=\ (\frac{2a}{\sqrt{3}},+\infty)^3$\Francois{, the lower bound being chosen in order to avoid overlaps of the balls,} $\mathcal{P} = \R^2\times \R $, and the functions $\bfX_i$ are now defined as
\Be
\bfX_i(\xi,\bfc,\alpha,\bfr) = \bfc + \Rot_\theta(\xi_i \bft_i + \bfr)\,\ \forall i \in\{1,2,3\}\,.
\Ee 
Notice that the functions $\bfX_i$ are 
still analytic in $(\xi,\bfc,\theta)$, and we use them to compute the
instantaneous velocity on the sphere $B_i$
\Be
\bfv_i=\frac{\partial \bfX_i}{\partial t}(\xi,\bfc,\theta,r) = \dot{\bfc}+\dot{\theta}\bfe_3 \times(\xi_i \bft_i + \bfr)+\Rot_\theta \bft_i \dot{\xi}_i\,,
\label{velocities}
\Ee
where $\bfe_3$  is the vertical unit vector. Eventually, due to the symmetries of the system, the swimmer stays in the
horizontal plane. 
\subsection{The four sphere swimmer moving in space (4S)}
\label{swimmer3}
We now turn to the more difficult situation of a swimmer able to move in the whole three dimensional space and rotate in any direction. In this case, we fix $N=4$ and we consider a regular reference
tetrahedron $(S_1,S_2,S_3,S_4)$ with center $O\in \R^3$ such
that $\mbox{dist}(O,S_i)=1$ and as before, we call $\bft_i = \vec{OS}_i$ for $i=1,2,3,4$.

The position and orientation in the three dimensional space of the tetrahedron 
are described by the coordinates of the
center $\bfc\in \R^3$ and a rotation $\mathcal{R}\in \SO$, in such a way that $d=6$.  

\begin{figure}[!htb]
\begin{center}
  \includegraphics[width=.35\textwidth]{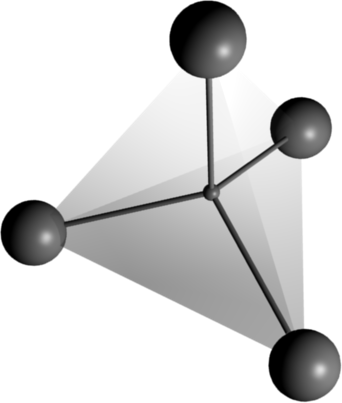}        
  \hfill
  
  \resizebox{.35\textwidth}{!}{
    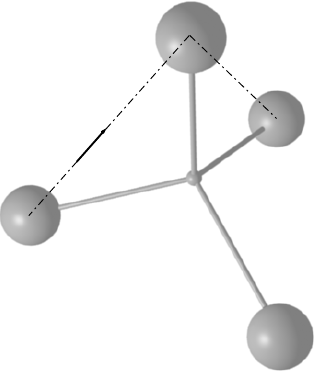
  }

  \caption{The four sphere swimmer (4S).}
  \label{fig:space-swimmer}
\end{center}
\end{figure}

We place the center of the ball $B_i$ at $\bfx_i=\bfc+\xi_i \Rot \bft_i$ with $\xi_i>0$ for $i=1,2,3,4$ as depicted in Fig. \ref{fig:space-swimmer} and forbid possible rotation of the spheres around the axes. A global rotation ($\Rot\ne \bfId$) of the swimmer is however allowed.

The four ball cluster is now completely described by the list of parameters
$X=(\xi,\bfc,\Rot)\in \mathcal{S}\times \mathcal{P}$, where $\mathcal{S} \ :=\ (\sqrt{\frac{3}{2}},+\infty)^4$ and $\mathcal{P} = \R^3\times \SO $. \Francois{Again, the lower bound for $\xi_i$ is chosen in order to avoid overlaps of the balls.}

Furthermore, the function $\bfX_i$ are now defined as
\Be
\bfX_i(\xi,\bfc,\Rot,\bfr) = \bfc + \Rot(\xi_i \bft_i + \bfr)\,\ \forall i \in\{1,2,3,4\}\,,
\Ee 
which \Francois{are} still analytic in $(\xi,\bfc,\Rot)$, from which we compute the
instantaneous velocity on the sphere $B_i$
\Be
\bfv_i=\frac{\partial \bfX_i}{\partial t}(\xi,\bfc,\Rot,r) = \dot{\bfc}+\bfomega\times(\xi_i \bft_i + \bfr)+\Rot \bft_i \dot{\xi}_i
\label{velocities-b}
\Ee
where $\bfomega$ is the angular velocity. Recall that $\bfomega=Ax(\dot{\Rot} \Rot)$ is the axial vector associated with the skew matrix $\dot{\Rot} \Rot$.

\Francois{
\subsection{The result}
In the theoretical part of this paper, we will establish the following theorem.
\begin{theorem}
Consider any of the swimmers described in Sections \ref{swimmer1}, \ref{swimmer2} or \ref{swimmer3}, and assume it is self-propelled in a three dimensional infinite viscous flow modeled by Stokes equations. Then for any initial configuration $(\xi^i,p^i)\in \mathcal{S}\times \mathcal{P}$ any final configuration $(\xi^f,p^f)\in \mathcal{S}\times \mathcal{P}$ and any final time $T>0$, there exists a stroke $\xi\in \mathcal{C}^0([0,T],\mathcal{S})$, piecewise $\mathcal{C}^1([0,T],\mathcal{S})$  satisfying
$\xi(0)=\xi^i$ and $\xi(T)=\xi^f$ such that if the self-propelled swimmer starts in position $p^i$ with the shape $\xi^i$ at time $t=0$, it ends at position $p^f$ and shape $\xi^f$ at time $t=T$ by changing its shape along $\xi(t)$.
\label{ourtheorem}
\end{theorem}
}

\section{Modelization of the fluid}
\label{sec2}

\Francois{In this section we give the expression of the total force and torque induced by
a shape change of the swimmer. Since the swimmer is composed of unions of balls, this turns out to study the Dirichlet-to-Neumann map outside the swimmer.} 
\Francois{We recall that  the swimmer is composed of $N$ identical and non-intersecting balls $B_1,\cdots,B_N$ of radius $a>0$ not necessarily aligned,
of center $\bfx_i\in \Re^3$, but linked together by deformable jacks which form a kind of skeleton. Since the lateral size of the arms of this skeleton are negligible, we can consider that the fluid fills the unbounded domain}
$\Omega=\Re^3\setminus \cup_{i=1}^N \overline{B}_i$. 
Since the balls do not intersect, the vector
$\bfx=(\bfx_1,\cdots,\bfx_N)\in \Omega$ is restricted to belong to
\Be
\AA :=\ \left\{ \bfx\in (\R^3)^N
\ :\ \min\limits_{i<j} |\bfx_i-\bfx_j| >2a\right\}\,.
\Ee

We work at low Reynolds number,  so that the fluid obeys Stokes equations
\Ben
\label{Stokes-Forces}
\left\{ \begin{array}{rcll}
-\eta \Delta \bfu + \nab p &=& 0& \mbox{ in }\Omega,\\
\nabdot \bfu &=& 0&\mbox{ in }\Omega,\\
-\bfsig \bfn &=& \bff &\mbox{ on }\pO,\\
\bfu&\to&0&\mbox{ at }\infty\,.
\end{array}\right.  
\Een
where $(\bfu,p)$ are respectively the velocity and the pressure of the fluid, $\eta$ its viscosity, $\bfsig := \eta(\nabla \bfu+ \nabla \bfu^t)- p \bfId$ is the Cauchy stress tensor and $\bfn$ is the outer unit normal to $\pO$
(hence, $\bfn$ points from the fluid to the interior of the balls). Existence and uniqueness of a solution to (\ref{Stokes-Forces}) is classical in the Hilbert space
\Be
\mathcal{V}\ :=\ \left\{ \bfu \in \mathcal{D}'(\Omega,\R^3)\ |
\ \nab \bfu \in L^2(\Omega) , \,\cfrac{\bfu}{\sqrt{1+|\bfr|^2}}\in L^2(\Omega)\right\},
\Ee
endowed with the norm 
\Be
\| \bfu\|_\mathcal{V}^2 \ :=\ \int_{\Omega} \left|\nabla \bfu\right|^2.
\Ee
Assuming that the force field belongs to
$\mathcal{H}^{-1/2}(\pO)$\footnote{Here, and in all the paper the symbol $\mathcal{H}^s$ denotes the Sobolev space of order $s$.}, the solution of (\ref{Stokes-Forces}) can be expressed in terms of the associated Green's function, namely the stokeslet
\begin{equation}
\bfG(\bfr)\ :=\ \cfrac1{8\pi\eta} \left( \cfrac1{|\bfr|} + \cfrac{\bfr\otimes \bfr}{|\bfr|^3}\right)
\label{stokeslet}
\end{equation}
as
\begin{equation}
\bfu(\bfr)=\int_{\pO}\bfG(\bfr-\sigma)\bff(\sigma)\,d\sigma
\label{fgivesu}
\end{equation}
where $\bff$ is a distribution of forces on $\partial \Omega$.
The Neumann-to-Dirichlet map
\begin{equation}
\begin{array}{rcl}
\mathcal{T}:\mathcal{H}^{-1/2}(\pO)&\longrightarrow& \mathcal{H}^{1/2}(\pO)\\
\bff & \longmapsto &\displaystyle  \bfu_{|\pO}=\left(\int_{\pO}\bfG(\bfr-\sigma)\bff(\sigma)\,d\sigma\right)_{|\pO}
\end{array}
\end{equation}
is a one to one mapping onto while, by the
open mapping Theorem, its inverse (the Dirichlet-to Neumann map) 
$\TT^{-1}$ is continuous. 

Calling $B=B(O,a)$ the ball of radius $a$ centered in the origin $O$, the special structure of our domain
allows us to identify $\pO$ with $(\pB)^N$ in such a way that $\TT$ can be seen as a map
\begin{equation}
\begin{array}{rcl}
\TT_\bfx:\mathcal{H}^{-1/2}_N&\longrightarrow& \mathcal{H}^{1/2}_N\\
(\bff_1,\cdots,\bff_N) & \longmapsto & (\bfu_1,\cdots,\bfu_N)
\end{array}
\end{equation}
in which $\mathcal{H}^s_N$ stands for
$\left(\mathcal{H}^s(\pB)\right)^N$ and the dependence on $\bfx\in
\AA$ has been emphasized, defined by
\begin{equation}
\begin{array}{rcl}
\bfu_i(\bfr) &=&\displaystyle \sum_{j=1}^N\int_{\partial B} \bfG(\bfx_i-\bfx_j+\bfr-\sigma)\bff_j(\sigma)\,d\sigma\\
&=:&\displaystyle \sum_{j=1}^N
\lra{\bff_j,\bfG(\bfx_i-\bfx_j+\bfr-\cdot)}_{\pB}, \quad \bfr \in \partial B.
\end{array}
\label{Green2}
\end{equation}

\subsection{Analyticity of the Dirichlet-to-Neumann map}
\label{section1}
This section is devoted to the following lemma which will be essential
in the proof of
the controllability theorem. \Francois{This kind of result was already proved in \cite{AlougesDeSimoneLefebvre-2007-a} but with a much more complicated method. We therefore give another much simpler proof which applies to more general situations, e.g., the swimmers described in subsections \ref{swimmer2} and \ref{swimmer3}.} We denote by $\mathcal{L}(E,F)$ the Banach space of linear maps
from $E$ to $F$ endowed with its usual norm.

\begin{lemma}
\label{lemanalytic}
The mapping  $\bfx \longmapsto \TT_\bfx$  is analytic from $\AA$ into $\mathcal{L}(\HH^{-1/2},\HH^{1/2})$.
\end{lemma}

The notion of analyticity in a Banach space is classical: it means
that at all points $\bfx_0 \in \AA$, $\TT_\bfx$ is equal to its Taylor series which converges in $\mathcal{L}(\HH^{-1/2},\HH^{1/2})$
for all $\bfx$ in a suitable neighborhood of $\bfx_0$.

\begin{proof}[Proof of Lemma~\ref{lemanalytic}]~
  
Let $\bfx\in \AA$. From \eqref{Green2}, we have for all $\bff=(\bff_1,\cdots,\bff_N)\in \HH^{-1/2}$
and for all $i \in \{1,\cdots,N\}$
\Be
(\TT_{\bfx}\bff)_{i}(r) \ =\  \lra{\bff_i ,\bfG(\bfr-\cdot)}_{\pB} + \sum_{j\neq
  i}\lra{\bff_j ,\bfG(\bfx_i-\bfx_j+\bfr-\cdot)}_{\pB}.
\Ee
Since the first term does not depend on $\bfx$, we only need to prove
the analyticity of
\Be
\begin{array}{rcl}
  \Psi \ :\ \R^3\setminus 2 B &\longrightarrow &  \mathcal{L}(\H^{-1/2}(\partial B),\H^{1/2}(\partial B)),\\
\bfz &\longmapsto &  \bigg(\: \bff \mapsto  \bfg(\bfr):=
\lra{\bff,\bfG(\bfz+\bfr-\cdot)}_{\pB}\: \bigg)
\end{array}
\Ee
which follows from the analyticity of $\bfG$ in
$\R^3\setminus\{0\}$ since this formula does not 
involve the singularity of the Green kernel
$\bfG$. 
\end{proof}
We have the following consequence.  
\Bc
Since  $\TT_\bfx$ is an isomorphism for every $\bfx\in \AA$, the mapping
$\bfx\longmapsto \TT_\bfx^{-1}$ is also analytic from $\AA$ to $\mathcal{L}(\HH^{1/2},\HH^{-1/2})$.  
\Ec
\begin{proof}
For $\bfx\in \AA$ and for $|\bfy|$ small enough, writing the Taylor series of $\TT_{\bfx+\bfy}$
as
\Be
\TT_{\bfx+\bfy}=\TT_\bfx+\sum_{q\geq 1} \TT_{\bfx,q}(\bfy)\,,
\Ee
where $\TT_{\bfx,q}(\bfy)$ denotes the $q-$homogeneous term in $\bfy$ of $\TT_{\bfx+\bfy}$, we have
\Be
 \TT_{\bfx+\bfy}^{-1} = \left\{ \left(Id + \sum_{q\geq 1}
 \TT_{\bfx,q}(y)\TT_\bfx^{-1}\right) \TT_\bfx \right\}^{-1} =
\TT_\bfx^{-1} \left\{ \sum_{p\geq 0} (-1)^p\left(\sum_{q\geq 1}  \TT_{\bfx,q}(y)\TT_\bfx^{-1}\right)^p\right\}
\Ee
where the last expression may be rearranged as a converging Taylor series.
\end{proof}

\subsection{Approximation for large distances}

We will also use the asymptotic behavior of $\TT_\bfx^{-1}$ as $\delta(\bfx):= \min_{i<j} |\bfx_i-\bfx_j|$ tends to $\infty$. 
\noindent
For $i,j \in \{1,\cdots,N\}$, we introduce the notation 
\Be
\bfx_{i,j}\ :=\ \bfx_i -\ \bfx_j,\qquad  \rij \ :=\  |\bfx_{i,j}| \qquad \mbox{and} \qquad \bfe_{i,j}\ := \ \cfrac{\bfx_{i,j}}\rij\,,
\Ee
and, for every $\bff\in \H^{-1/2}(\partial B,\R^3)$,
\Be
\left(\overline{\TT} \bff\right)(\bfr) = \int_{\partial B}\bfG(\bfr-\sigma)\bff(\sigma)\,d\sigma\,, \quad \forall \bfr\in \partial B\,.
\Ee
\begin{proposition}
\label{prop_approx_Phi}
For every $\bff\in \HH^{-1/2}$, 
\Ben
(\TT_\bfx \bff)_{i} \ =\ \overline{\TT}  \bff_i + \cfrac1{8\pi\eta}\sum_{j\neq i}
\cfrac1\rij\left( \bfId+ \bfe_{i,j}\otimes\bfe_{i,j}\right) \lra{\bff_j,\bfId}_{\pB} + O(\delta^{-2})||\bff||_{\HH^{-1/2}}\,.  
\label{expT}
\Een
For every $\bfg\in \HH^{1/2}$, $\bfx\in \AA$ and $\bfy\in (\R^3)^N$, we have 
\Ben
(\TT_\bfx^{-1} \bfg)_{i} \ =\ \overline{\TT}^{-1}  \left\{ \bfg_i - \cfrac1{8\pi\eta}\sum_{j\neq i}
\cfrac1\rij\left( \bfId+ \bfe_{i,j}\otimes\bfe_{i,j}\right) \lra{ \overline{\TT}^{-1} \bfg_j
  ,\bfId}_{\pB}\right\} + O(\delta^{-2})||\bfg||_{\HH^{1/2}}\,. 
\label{expTminus1}
\Een
\end{proposition}
\begin{proof}
Let $\bfx\in \AA$
and $\bff\in \HH^{-1/2}$. From~\eqref{Green2} we
have for $i=1,\cdots,N$ and $\bfr\in \pB$,  
\Ben
(\TT_\bfx \bff)_{i}(\bfr) \ =\ \overline{\TT} \bff_i (\bfr)  + \sum_{j\neq i} \lra{\bff_j,\bfG(\bfx_{i,j}+\bfr-\cdot)}_{\pB}\,.
\label{expansion}
\Een
Using the expansion 
\Be
\bfG(\bfx_{i,j}+ \bfz) = \bfG(\bfx_{i,j}) + O(\delta^{-2})\,,
\Ee
which is valid  uniformly for bounded $\bfz$,
we get
\Be
\lra{\bff_j,\bfG(\bfx_{i,j}+\bfr-\cdot)}_{\pB}\ =\ \bfG(\bfx_{i,j})\, \lra{\bff_j,\bfId}_{\pB} + O(\delta^{-2})||\bff||_{\HH^{-1/2}}\,.
\Ee
Substituting into (\ref{expansion}), and using (\ref{stokeslet}), we obtain (\ref{expT}).
For (\ref{expTminus1}) we remark that (\ref{expT}) is the asymptotic expansion in $\bfx$ of the analytic operator $\TT_\bfx$ at infinity. Since $\overline{\TT}$ is invertible, we obtain the expansion (\ref{expTminus1}) by inverting (\ref{expT}).
\end{proof}

In the sequel, we will use furthermore the fact that the spheres are
non-deformable and may thus only move following a rigid body motion. In
this case, the velocity of each point $\bfr$ of the $i-$th sphere is given by
\begin{equation}
\bfv_i(\bfr) = \bfu_i+\omega_i\times \bfr\,.
\label{localvelocity}
\end{equation}
Using (\ref{expTminus1}) we find the expansion of the total force that the system applies to the fluid
\Ben
\bfF = \sum_{i=1}^N 6\pi\eta a\bfu_i -  \frac{9\pi\eta a^2}{2}\sum_{i=1}^N\sum_{j\neq i}
\bfS_{ij}\bfu_j+ O(\delta^{-2})||\bfu||\,,
\label{totalforce}
\Een
where $\bfS_{ij}$ is the matrix defined by
\begin{equation}
\bfS_{ij} = \cfrac1\rij\left( \bfId+ \bfe_{i,j}\otimes\bfe_{i,j}\right)\,,
\label{sij}
\end{equation}
while the total torque (with respect to the origin $O$) is given by
\Ben
\bfT = 6\pi\eta a\sum_{i=1}^N \bfx_i\times \bfu_i + O(1)||\bfu||\,.
\label{totaltorque}
\Een

\Francois{
\section{The swimming problem as an affine control problem without drift}
\label{sec3}
In this section we show that self-propulsion of the swimmers makes the system that
describes their dynamics a linear control system.
\begin{lemma}The system that describes the dynamics of low $Re$ swimmers can be written in the form
\begin{equation}
\frac{d}{dt} \left(\begin{array}{c} \xi\\p\end{array}\right) = \sum_{i=1}^M \bfF_i(\xi,p)\dot{\xi}_i\,,
\label{controleq}
\end{equation}
where the vectorfields $(\bfF_i)_{1\leq i\leq M}$ are analytic in $(\xi,p)$. It is therefore an affine control problem without drift where the controls are the rate of shape changes.
\end{lemma}}

\Francois{The rest of this section is devoted to the proof of this lemma. This can be seen as a generalization of the result obtained in \cite{AlougesDeSimoneLefebvre-2007-a} for the three sphere swimmer of Najafi and Golestanian. We show in particular how self-propulsion gives rise to the specific form (\ref{controleq}) of the equation, describing the dynamics of all these swimmers of Section \ref{sec1}.}

\subsection{Self-propulsion}

The question we want to address is whether it is possible to control
the state of the system (i.e. both $\xi$ and $p$) using as controls
only the rate of shape changes $\dot{\xi}$. For this aim, we need to
understand the way $p$ varies when one changes $\xi$. This is done
assuming that the swimmer is self-propelled, and that the swimmer's inertia
is negligible, which implies that the total viscous force and torque
exerted by the surrounding fluid on the swimmer must vanish.  As we
shall see, using these conditions, $\dot{p}$ is uniquely \Francois{and linearly} determined by
$\dot{\xi}$. The condition that total viscous force and torque vanish
is written as \Ben \sum_{i=1}^N \int_{\pB} \TT^{-1}_{\bfx(\xi,p)}
\bfv_i(\xi,p,\bfr)\,d\sigma(\bfr) = 0\,,
\label{noforce}
\Een
and 
\Ben
\sum_{i=1}^N \int_{\pB} \bfX_i(\xi,p,\bfr)\times \TT^{-1}_{\bfx(\xi,p)} \bfv_i(\xi,p,\bfr)\,d\sigma(\bfr) = 0\,,
\label{notorque}
\Een
which lead, using (\ref{eqvi}), to a set of linear equations which
link $\dot{p}$ to $\dot{\xi}$. In all the examples we have in mind
(see below),
this allows us to solve $\dot{p}$ uniquely and linearly in terms of $\dot{\xi}$
\begin{equation}
\label{linear-relationship-xi-dot-p}
\dot{p} = \sum_{i=1}^M \VW_i(\xi,p)\dot{\xi}_i\,.
\end{equation}
In the preceding equation $(\VV_i(\xi,p))_{1\leq i\leq M}$ are $M$ vector fields 
defined on $T\mathcal{P}$, the tangent bundle of $\mathcal{P}$.

The state of the system $(\xi,p)$ thus follows the ODE
\begin{equation}
\left(
\begin{array}{c} \dot{\xi}\\\dot{p}\end{array}\right) = \sum_{i=1}^M \FF_i(\xi,p)\dot{\xi}_i\,,
\label{controlsystems}
\end{equation}
where now for all $1\leq i\leq M$, $\FF_i(\xi,p) =
\left(\begin{array}{c} \bfe_i\\  \VW_i(\xi,p)\end{array}\right)$
is a vector field defined on $T\mathcal{S}\times T\mathcal{P}$ (here $\bfe_i$ is the $i$-th vector of
the canonical basis of $\Re^M$).

\subsection{The three sphere swimmer of Najafi and Golestanian (3S)}

Due to axial symmetry, the only non trivial equation in (\ref{noforce}) is the one
concerning the component of the total viscous force along the axis of symmetry.
This reads
\Be
\sum_{k=1}^2 \lambda_k(\xi,p) \dot{\xi}_k+\lambda_0(\xi,p) \dot{p} =0,
\Ee
where
\Be
\begin{array}{rcl}
\lambda_1(\xi) &=&\displaystyle \sum_{i=1}^3\int_{\pB} \bfe_1\cdot \left(\TT^{-1}_{\bfx(\xi,p)}\left(\frac{\partial \bfX_1}{\partial \xi_1},\frac{\partial \bfX_2}{\partial \xi_1},\frac{\partial \bfX_3}{\partial \xi_1}\right)\right)_i\,d\bfr\\
&=&\displaystyle \sum_{i=1}^3\int_{\pB} \frac{\partial \bfX_i}{\partial \xi_1}\cdot \left(\TT^{-1}_{\bfx(\xi,p)}\left(\bfe_1,\bfe_1,\bfe_1\right)\right)_i\,d\bfr\\
&=&\displaystyle -\int_{\pB} \bfe_1\cdot \left(\TT^{-1}_{\bfx(\xi,p)}(\bfe_1,\bfe_1,\bfe_1)\right)_1\,d\bfr\,,
\end{array}
\Ee
and, similarly,
\Be
\lambda_2(\xi) = \displaystyle \int_{\pB} \bfe_1\cdot \left(\TT^{-1}_{\bfx(\xi,p)}(\bfe_1,\bfe_1,\bfe_1)\right)_3\,d\bfr\,,
\Ee
while the coefficient of $\dot{p}$ is simply given by
\Be
\begin{array}{rcl}
\lambda_0(\xi) &=& \displaystyle \sum_{i=1}^3\int_{\pB} \bfe_1\cdot \left(\TT^{-1}_{\bfx(\xi,p)}(\bfe_1,\bfe_1,\bfe_1)\right)_i\,d\bfr\\
&>& 0\,.
\end{array}
\Ee
Positivity of the viscous drag $\lambda_0$ arising from a rigid translation at unit speed is a consequence of the positive-definiteness of the Oseen resistance matrix, see \cite{HappelBrenner-1965-a}.
 The coefficients above are independent of $p$ because of translational invariance of the problem. From the analyticity of the coefficients and $\TT_\bfx$
in $\bfx$, we have that 
\Be
\dot{p}=\sum_{k=1}^2 V_{k}(\xi)\dot{\xi}_k\,,
\Ee
where $\displaystyle V_k(\xi)=-\frac{\lambda_k(\xi)}{\lambda_0(\xi)}$
are analytic functions of $\xi_1,\xi_2$. The complete system reads
\Ben
\frac{d}{dt}\left(\begin{array}{c}\xi\\p\end{array}\right) = \alpha_1(t)
\left(\begin{array}{c}1\\0\\ V_1(\xi)\end{array}\right)+\alpha_2(t)
\left(\begin{array}{c}0\\1\\ V_2(\xi)\end{array}\right)
\label{system3S}
\Een
where $\alpha_1(t)=\dot{\xi}_1(t)$ and $\alpha_2(t)=\dot{\xi}_2(t)$ are the control functions. This is under the form (\ref{controleq}).

\subsection{The three sphere swimmer moving in a plane (3SP)}

Due to the symmetries of the system, the swimmer stays in the
horizontal plane. Therefore, the third component of the total force
$F_3$ vanishes, while only the third component of the total torque
$T_3$ might not vanish. The constraints of self-propulsion
(\ref{noforce}) and (\ref{notorque}) -- more precisely, the vanishing
of the components $F_1$, and $F_2$ of the force, and of the third
component $T_3$ of the torque -- lead to a linear system that can be
written as \Be \Res(\xi,\theta)\left(\begin{array}{c} \dot{\bfc}\\
    \dot{\theta}
\end{array}
\right) + \sum_{i=1}^3 \VV_i(\xi,\theta) \dot{\xi}_i = 0\,,
\Ee
where neither $\Res$, nor $\VV_i$ depend on $\bfc$ because of translational invariance of the Stokes problem.

The complete system reads
\Ben
\left(\begin{array}{c}\dot{\xi} \\ \dot{\bfc} \\ \dot{\theta}\end{array}\right) = 
\left(\begin{array}{c}\bfId\\ -\Res^{-1}(\xi,\theta)\VV(\xi,\theta)
\end{array}\right)\alpha(t)\,,
\label{system3SP}
\Een
where $\VV(\xi,\theta)$ is the $3\times 3$ matrix whose columns are $\VV_i(\xi,\theta)$,
and $\alpha(t)=\dot{\xi}\in \R^3$ are the controls. This is again under the form (\ref{controleq}).

\subsection{The four sphere swimmer moving in space (4S)}
Constraints of self propulsion (\ref{noforce}) and (\ref{notorque}) now become
\Be
\left(\begin{array}{cc}\bfK(\xi,\Rot) & \bfC^T(\xi,\Rot)\\
\bfC(\xi,\Rot) & \bfOm(\xi,\Rot)
\end{array}
\right)
\left(\begin{array}{c}
\dot{\bfc}\\ \bfomega
\end{array}
\right) + \sum_{i=1}^4 \VV_i(\xi,\Rot) \dot{\xi}_i = 0\,,
\label{systemtotal}
\Ee
where the $6\times 6$ matrix
\Ben
\Res(\xi,\Rot) = \left(\begin{array}{cc}\bfK(\xi,\Rot) & \bfC^T(\xi,\Rot)\\
\bfC(\xi,\Rot) & \bfOm(\xi,\Rot)
\end{array}
\right)
\label{resistance}
\Een
known as the resistance matrix \cite{HappelBrenner-1965-a}
is symmetric positive definite,
and depends analytically on $(\xi,\Rot)$, just as the vectors $\VV_i(\xi,\Rot)$ do.

In this case, the complete system reads
\Ben
\left(\begin{array}{c}\dot{\xi} \\ \dot{\bfc} \\ \bfomega\end{array}\right) = 
\left(\begin{array}{c}\bfId_{4,4}\\ -\Res^{-1}(\xi,\Rot)\VV(\xi,\Rot)
\end{array}\right)\alpha(t)\,,
\label{control4S1}
\Een
where $\VV(\xi,\Rot)$ is the $6\times 4$ matrix whose columns are $\VV_i(\xi,\Rot)$,
and $\alpha(t)\in \R^4$ are the controls. Using $\dot{\mathcal{R}}=\mbox{skew}(\bfomega)\mathcal{R}$ enables us to write the system under the more classical form (\ref{controleq})
\Ben
\left(\begin{array}{c}\dot{\xi} \\ \dot{\bfc} \\ \dot{\mathcal{R}}\end{array}\right) = 
\sum_{i=1}^4\bfF_i(\xi,\mathcal{R})\alpha_i(t)\,,
\label{control4S2}
\Een
where $\alpha_i(t)=\dot{\xi}_i(t)$ are the control functions.

\Francois{\section{Proof of theorem \ref{ourtheorem}}
\label{sec4}
}
\subsection{Geometric control theory}

When we use $(\dot{\xi}_i)_{1\leq i \leq M}$ as control variables, systems like (\ref{controlsystems}) belong to the class of linear control systems with analytic
coefficients of the form
\Ben
\frac{dX}{dt}=\sum_{i=1}^M \alpha_i(t) F_i(X)
\label{generalsystem}
\Een
where $X$ is a point of a $d-$dimensional manifold $\mathcal{M}$, and $(F_i)_{1\leq i\leq M}$ are $p$ vector fields defined on $T\mathcal{M}$, the tangent bundle of $\mathcal{M}$. For such systems, the tools of geometric control theory can be applied. In our case, $\mathcal{M}$ and the vector fields $F_i$ are furthermore analytic, which enables us to use stronger results. The theory
for such systems may be found in \cite{AgrachevSachkov-2004-a,Jurdjevic-1997-a}, for instance.

In (\ref{generalsystem}), the control which governs the system is the vector $(\alpha_1,\cdots,\alpha_p)$, and the question of controllability of the system can be stated as follows:

\begin{question} For any \Francois{pair} of states $(X_0,X_1)\in \mathcal{M}^2$, are there $M$ functions $\alpha_i:[0,T]\rightarrow \R$ such that solving
(\ref{generalsystem}) starting from $X(0)=X_0$ leads to $X(T)=X_1$?
\label{Q1}
\end{question}

A classical result by Rashevsky and Chow states that the system (\ref{generalsystem}) is locally
controllable near $X_0\in \mathcal{M}$ (which means that question \ref{Q1} 
possesses a positive answer for any $X_1$ in a suitable neighborhood of $X_0$) \Francois{ with piecewise constant controls $\alpha_i$ that are moreover continuous from the right}, 
if the Lie algebra generated by $(F_1,\cdots,F_M)$ is of full rank at $X_0$
\Ben
\mbox{dim }Lie(F_1,\cdots,F_M)(X_0)= d\,.
\label{fullrank}
\Een
Here, the Lie algebra $Lie(F_1,\cdots,F_M)(X_0)$ is the algebra obtained
from the vector fields $(F_1,\cdots,F_M)$ by successive applications of the Lie
bracket operation defined as
\Be
[F,G](X)=(F\cdot \nabla) G(X) - (G\cdot \nabla) F (X)\,.
\Ee
The rigorous proof of Rashevsky-Chow's theorem may be found in \cite{AgrachevSachkov-2004-a,Jurdjevic-1997-a} but we 
would like to give here a hint of why such a result holds.
Given a vector field $F:\mathcal{M} \longrightarrow T\mathcal{M}$, we call $\exp(tF)(X_0)$ 
the solution at time $t$ of
\Ben
\left[
\begin{array}{l}
\displaystyle \frac{dX}{dt}=F(X)\,,\\
X(0)=X_0\,.
\end{array}
\right.
\label{expfunction}
\Een
The main idea is that one can reach from $X_0$ points in the direction
$g = \sum_i \beta_i F_i(X_0)$ by using (\ref{expfunction}) with $F=\sum_i \beta_i F_i$ since we have, to first order, for small $|t|<<1$,
\Be
\exp(t F)(X_0)=X_0+t F(X_0)+O(t^2)=X_0+t g + O(t^2)\,.
\Ee
This formula shows that the global displacement of the solution of the dynamical system is proportional to time for these directions.

A more subtle move allows \Francois{one} to reach points in the direction $[F_i,F_j]$.
Namely, we now need to consider
\Be
\exp(-tF_i)\circ \exp(-tF_j)\circ \exp(tF_i)\circ \exp(tF_j)(X_0) = X_0 + t^2[F_i,F_j](X_0) + O(t^3)\,,
\Ee 
which also shows that in time $t$, one can reach a displacement in the desired direction which is proportional to $t^2$. Lie bracket directions 
of higher order are also attainable at the price of higher order expressions 
in $t$, leading to smaller displacements. If the Lie algebra has a full rank,
every direction in $\mathcal{M}$ is attainable and the system is locally controllable.

We will focus on systems 
for which the full rank condition (\ref{fullrank}) is satisfied
with first order Lie brackets only, and a necessary condition for
\Be
\mbox{dim Span}(F_1,\cdots,F_M,[F_1,F_2],[F_1,F_3],\cdots,[F_{M-1},F_M])(X_0)=d\,,
\Ee
is that 
\Be
M+\frac{M(M-1)}{2}\geq d\,.
\Ee
Our three systems are designed in such a way that they all
satisfy $\displaystyle M+\frac{M(M-1)}{2}= d$, and  local
controllability near $X_0$ follows from verifying that
\Ben
\mbox{det }(F_1,\cdots,F_M,[F_1,F_2],[F_1,F_3],\cdots,[F_{M-1},F_M])(X_0)\ne 0\,.
\label{sufcond}
\Een

The remaining
of this section consists in the detailed analysis of \Francois{our three model swimmers. In each case, we show that (\ref{sufcond}) holds at one point.}
proving local controllability at one point
$(\xi_0,p_0)$. In order to obtain global controllability\Francois{, and theorem \ref{ourtheorem}, we need a special argument which is the aim of Section \ref{localtoglobal}. The regularity of the shape function $\xi(t)$ follows from the fact that the controls are here the rate of shape changes.}  


\subsection{The 3S swimmer}
For the equation (\ref{system3S}), the controllability condition (\ref{sufcond}) becomes
\begin{equation}
\mbox{det}(F_1,F_2,[F_1,F_2])(\xi)=\frac{\partial V_2}{\partial \xi_1}(\xi)-\frac{\partial V_1}{\partial \xi_2}(\xi) \neq 0\,,
\label{sufcond1}
\end{equation}
which is enough to verify for one configuration. For this purpose, and
in order to further simplify the computation, we make the observation
that, due to symmetry, 
\Be
\lambda_1(\xi_1,\xi_2)=-\lambda_2(\xi_2,\xi_1),\mbox{ while
}\lambda_0(\xi_1,\xi_2)=\lambda_0(\xi_2,\xi_1)\,.  \Ee Therefore, if
we consider a symmetric shape $\overline{\xi}=(\zeta,\zeta)$, showing
that condition (\ref{sufcond1}) holds at $\overline{\xi}$ is
equivalent to proving that \Be \frac{\partial \lambda_1}{\partial
  \xi_2}(\overline{\xi})\neq 0\,.  \Ee Using the expansion for large
distances of the total force (\ref{totalforce}) with
$\bfu_1=(\dot{p}-\dot{\xi}_1) \bfe_1$, $\bfu_2=\dot{p} \bfe_1$, and
$\bfu_3=(\dot{p}+\dot{\xi}_2) \bfe_1$, leads to
\begin{eqnarray*}
&&\lambda_0(\xi_1,\xi_2)=6\pi a \eta\left(3-\frac{3a}{2}\left(\frac{1}{\xi_1}+\frac{1}{\xi_2}+\frac{1}{\xi_1+\xi_2}\right)\right) +O(\zeta^{-2})\,, \\
&&\lambda_1(\xi_1,\xi_2)=6\pi a \eta\left(1-\frac{3a}{2}\left(\frac{1}{\xi_1}+\frac{1}{\xi_1+\xi_2}\right)\right) +O(\zeta^{-2})\,, \\
&&\lambda_2(\xi_1,\xi_2)=6\pi a \eta\left(-1+\frac{3a}{2}\left(\frac{1}{\xi_2}+\frac{1}{\xi_1+\xi_2}\right)\right)+O(\zeta^{-2})\,. 
\end{eqnarray*}
Therefore 
$$
\frac{\partial \lambda_1}{\partial \xi_2}(\overline{\xi}) = \frac{9\pi a^2 \eta}{4\zeta^2}+O(\zeta^{-3})
$$
which does not vanish for $\zeta$ sufficiently large. This proves the global controllability of the system.

\subsection{The 3SP swimmer}
We check the controllability condition (\ref{sufcond}) for the equation (\ref{system3SP}) at $\bar{\theta}=0$ and at a symmetric shape
$\bar{\xi}=(\zeta,\zeta,\zeta)$, with $\zeta \in \R$ sufficiently large. 

Calling $\bfF_i=-\Res^{-1}(\xi,\theta)\VV_i(\xi,\theta)$, the special form of our system (in particular the fact that the coefficients do not depend on the position variable $\bfc$) allows us to rewrite the sufficient condition (\ref{sufcond}) for
controllability at $(\bar{\xi},\bfc,\theta)$ as
\begin{equation}
\Delta = \mbox{det}\left(\frac{\partial \bfF_1}{\partial \xi_2}-\frac{\partial \bfF_2}{\partial \xi_1}, \frac{\partial \bfF_2}{\partial \xi_3}-\frac{\partial \bfF_3}{\partial \xi_2},\frac{\partial \bfF_3}{\partial \xi_1}-\frac{\partial \bfF_1}{\partial \xi_3}\right)(\bar{\xi})\ne 0.
\label{deter}
\end{equation}
Notice that the absence of $\partial/\partial \theta$ terms come from the fact that $(\bfF_i)_3=0$ in a symmetric configuration.
We now turn to the computation of $\displaystyle \frac{\partial \bfF_k}{\partial \xi_l}(\bar{\xi})$ for $k\ne l$ or, at least, of its asymptotic expansion in terms of negative powers of $\zeta$.

From the definition of $\bfF_k$, one has
\begin{equation}
\frac{\partial \bfF_k}{\partial \xi_l} = -\Res^{-1}(\xi,\theta)\left(\frac{\partial \VV_k}{\partial \xi_l}+\frac{\partial \Res}{\partial \xi_l}\bfF_k \right)\,,
\label{dfdxi}
\end{equation} 
and we need to compute all the factors of this expression at $(\bar{\xi},0)$.

The first step consists in computing an expansion for
$\Res(\bar{\xi},0)$. We get from (\ref{totalforce}) and
(\ref{totaltorque}), and taking $N=3$,
\begin{equation}
\Res(\bar{\xi},0)=6\pi a \eta\left(
\begin{array}{cc}
 3 \bfId +O(\zeta^{-1})& O(1)\\
O(1) & \displaystyle 3\zeta^2 + O(\zeta)
\end{array}
\right)\,,
\end{equation}
and
\begin{equation}
\Res^{-1}(\bar{\xi},0)=\frac{1}{6\pi a \eta}\left(
\begin{array}{cc}
\displaystyle \frac{1}{3} \bfId +O(\zeta^{-1})& O(\zeta^{-2})\\
O(\zeta^{-2}) & \displaystyle \frac{1}{3\zeta^2}  + O(\zeta^{-3})
\end{array}
\right)\,.
\end{equation}
Differentiating $\Res$ with respect to $\xi_l$ gives (notice that from
Section \ref{section1}, $\Res$ is analytic in $\xi$ and we may
therefore differentiate its expansion  term by term)
\begin{equation}
\frac{\partial \Res}{\partial \xi_l}(\bar{\xi},0)=6\pi a \eta\left(
\begin{array}{cc}
-\frac{3}{2}a \frac{\partial}{\partial \xi_l}\left(\sum_{i\ne l} \bfS_{il}\right) + O(\zeta^{-3})& \,\,\,\,\bft_l\times \bfe_3   +O(\zeta^{-1})\\
(\bft_l\times \bfe_3)^t +O(\zeta^{-1}) & O(\zeta)
\end{array}
\right)\,,
\end{equation}
where $\bfS_{ij}$, originally given by (\ref{sij}), is now seen as its
restriction to the horizontal plane, and the vector $\bft_l\times \bfe_3$ as a horizontal vector.

Similarly
\begin{equation}
\VV_k(\bar{\xi},0)=\left(
\begin{array}{l}
6\pi\eta a\bft_k-9\pi\eta a^2\sum_{i\ne k} \bfS_{ik} \bft_k+O(\zeta^{-2})\\
O(1)
\end{array}
\right)
\end{equation}
gives
\begin{equation}
\bfF_k(\bar{\xi},0)=-\frac{1}{3}\left(
\begin{array}{l}
\bft_k+O(\zeta^{-1})\\
0
\end{array}
\right)
\end{equation}
(here we have used the fact that, for a symmetric configuration, the
displacement of one sphere induces no torque) and
\begin{equation}
\frac{\partial \VV_k}{\partial{\xi_l}}(\bar{\xi},0)=-9\pi\eta a^2\left(
\begin{array}{l}
\displaystyle \frac{\partial \bfS_{lk}}{\partial \xi_l} \bft_k+O(\zeta^{-3})\\
O(\zeta^{-1})
\end{array}
\right)\,.
\end{equation}

From these expressions, and using (\ref{dfdxi}), and the identity
$$
\frac{\partial \bfS_{il}}{\partial \xi_l} = \frac{1}{r_{il}^2}\left(-(\bfe_{il}\cdot \bft_l)\bfId+\bft_l\otimes \bfe_{il}+\bfe_{il}\otimes \bft_l - 3 (\bfe_{il}\cdot \bft_l)\bfe_{il}\otimes \bfe_{il}\right)
$$
one gets (after some calculations) 
\begin{equation}
\frac{\partial \bfF_k}{\partial \xi_l} (\bar{\xi},0)-\frac{\partial \bfF_l}{\partial \xi_k} (\bar{\xi},0) = \left(
\begin{array}{l}
\frac{ 4a}{\zeta^2}(\bft_l-\bft_k)\\
-\frac{1}{9\zeta^2}\bft_l\times \bft_k\cdot \bfe_3
\end{array}
\right)+O(\zeta^{-3})\,.
\end{equation}
The determinant $\Delta$ given by (\ref{deter}) is thus equal to
\Be
\Delta = -\zeta^{-6}\frac{16 a^2}{9} \det\left(\begin{array}{cccc}
\bft_2 -\bft_1 & \bft_3-\bft_2& \bft_1-\bft_3\\
\bft_2\times \bft_1\cdot \bfe_3 & \bft_3\times \bft_2\cdot \bfe_3 & \bft_1\times \bft_3\cdot \bfe_3
\end{array}\right)+O(\zeta^{-7})
\Ee
which does not vanish for $\zeta$ large enough.

We emphasize again that the global controllability result proven here
means that, by using suitable controls, the swimmer is capable of
moving anywhere in the horizontal plane and of rotating around the vertical
axis by any angle.

\subsection{The 4S swimmer}
As before, we check the controllability condition (\ref{sufcond}) for the equation (\ref{control4S2}) at a symmetric shape
$\bar{\xi}=(\zeta,\zeta,\zeta,\zeta)$, with $\zeta$
sufficiently large. Also, we use $\Rot=\bfId$ in the
verification of the Lie bracket condition~\eqref{sufcond}.

Notice that for the symmetric configuration $\bar \xi $,
the vectors $\bfF_i$ defined by (\ref{control4S1},\ref{control4S2}), have no torque components, and that $\bfF_i$ does not depend on
$\bfc$ because of the translational invariance of the problem. Taking into account these facts, the sufficient condition  (\ref{sufcond}) for
controllability at $(\bar{\xi},0,\bfId)$ simplifies to
\begin{equation}
\mbox{dim Span}\left(\frac{\partial \bfF_1}{\partial \xi_2}-\frac{\partial \bfF_2}{\partial \xi_1}, \cdots,\frac{\partial \bfF_3}{\partial \xi_4}-\frac{\partial \bfF_4}{\partial \xi_3}\right)(\bar{\xi},\bfId)= 6.
\label{deter20}
\end{equation}

Notice that in the preceding equation, each vector is a tangent vector to the group of 3 dimensional rigid transformations, $SE(3)=\mathbb{R}^3\times \SO$. In order to simplify the computation, it is convenient to call $\bfphi_i=-\Res^{-1}(\bar{\xi},\Rot)\VV_i$ (see (\ref{control4S1})), and work in {\em translations and angular velocities}.
A simple change of coordinates shows that condition (\ref{deter20}) simply rewrites
\begin{equation}
\Delta = \mbox{det}\left(\frac{\partial \bfphi_1}{\partial \xi_2}-\frac{\partial \bfphi_2}{\partial \xi_1}, \cdots,\frac{\partial \bfphi_3}{\partial \xi_4}-\frac{\partial \bfphi_4}{\partial \xi_3}\right)(\bar{\xi},\bfId)\ne 0.
\label{deter2}
\end{equation}
We now turn to the computation of $\displaystyle \frac{\partial
  \bfphi_k}{\partial \xi_l}(\bar{\xi},\bfId)$ for $k\ne l$, or at
least of  its asymptotic expansion in terms of negative powers of $\zeta$.

The definition of $\bfphi_k$ still gives
\begin{equation}
\frac{\partial \bfphi_k}{\partial \xi_l} = -\Res^{-1}(\bar{\xi},\Rot)\left(\frac{\partial \VV_k}{\partial \xi_l}+\frac{\partial \Res}{\partial \xi_l}\bfphi_k \right)\,,
\label{dfdxi-b}
\end{equation} 
and we need to compute all the factors of this expression at $(\bar{\xi},\bfId)$.

The first step consists in computing an expansion of $\Res(\bar{\xi},\bfId)$. We have from (\ref{totalforce}) and (\ref{totaltorque})
\begin{eqnarray*}
\bfK &=&\displaystyle 6\pi a\eta N \bfId -  \frac{9\pi\eta a^2}{2}\sum_{i=1}^N\sum_{j\ne i}\bfS_{ij}+O(\zeta^{-2})\,,\\
\bfC &=&6\pi a \eta\sum_{i=1}^N \xi_i \mbox{ skew}(\bft_i) +O(\zeta^{0})\,,\mbox{ with } \mbox{skew}(\bft)\bfx = \bft \times \bfx \,,\,\,\forall \bfx \in \R^3\,,\\
\bfOm &=& 6\pi a \eta\sum_{i=1}^N \xi_i^2(\bfId- \bft_i\otimes \bft_i)+O(\zeta)\,.
\end{eqnarray*}
This gives (here $N=4$, $\sum_{i=1}^4 \bft_i=0$ and  $\sum_{i=1}^4 (\bfId- \bft_i\otimes \bft_i)=\frac{8}{3}\bfId$)
\begin{equation}
\Res(\bar{\xi},\bfId)=6\pi a \eta\left(
\begin{array}{cc}
 4 \bfId +O(\zeta^{-1})& O(\zeta^{-1})\\
O(\zeta^{-1}) & \displaystyle \frac{8\zeta^2}{3} \bfId + O(\zeta)
\end{array}
\right)\,,
\end{equation}
and
\begin{equation}
\Res^{-1}(\bar{\xi},\bfId)=\frac{1}{6\pi a \eta}\left(
\begin{array}{cc}
 \frac{1}{4} \bfId +O(\zeta^{-1})& O(\zeta^{-2})\\
O(\zeta^{-2}) & \displaystyle \frac{3}{8\zeta^2} \bfId + O(\zeta^{-3})
\end{array}
\right)\,.
\end{equation}
Differentiating $\Res$ with respect to $\xi_l$ gives (notice that from
Section \ref{section1}, $\Res$ is analytic in $\xi$ and we may
therefore differentiate its expansion term by term)
\begin{equation}
\frac{\partial \Res}{\partial \xi_l}(\bar{\xi},\bfId)=6\pi a \eta\left(
\begin{array}{cc}
-\frac{3}{2}a \frac{\partial}{\partial \xi_l}\left(\sum_{i\ne l} \bfS_{il}\right) + O(\zeta^{-3})& \,\,\,\,-\mbox{skew}(\bft_l)  +O(\zeta^{-1})\\
\mbox{skew}(\bft_l) +O(\zeta^{-1}) & O(\zeta)
\end{array}
\right)\,.
\end{equation}
Similarly
\begin{equation}
\VV_k=\left(
\begin{array}{l}
6\pi\eta a\bft_k-9\pi\eta a^2\sum_{i\ne k} \bfS_{ik} \bft_k+O(\zeta^{-2})\\
O(1)
\end{array}
\right)
\end{equation}
leads to
\begin{equation}
\bfphi_k(\bar{\xi},\bfId)=-\frac{1}{4}\left(
\begin{array}{l}
\bft_k+O(\zeta^{-1})\\
0
\end{array}
\right)
\end{equation}
(here we have used the fact that for a symmetric configuration, the
displacement of one sphere induces no torque) and
\begin{equation}
\frac{\partial \VV_k}{\partial{\xi_l}}(\bar{\xi},\bfId)=-9\pi\eta a^2\left(
\begin{array}{l}
\displaystyle \frac{\partial \bfS_{lk}}{\partial \xi_l} \bft_k+O(\zeta^{-3})\\
O(\zeta^{-1})
\end{array}
\right)\,.
\end{equation}

From these expressions, and using (\ref{dfdxi-b}), one gets after a tedious but straightforward computation 
\begin{equation}
\frac{\partial \bfphi_k}{\partial \xi_l} (\bar{\xi},\bfId)-\frac{\partial \bfphi_l}{\partial \xi_k} (\bar{\xi},\bfId) = \left(
\begin{array}{l}
\frac{9\sqrt{3} a}{64\sqrt{2}\zeta^2}(\bft_l-\bft_k)\\
\frac{3}{16\zeta^2}\bft_l\times \bft_k
\end{array}
\right)+O(\zeta^{-3})\,.
\end{equation}
The determinant $\Delta$ given by (\ref{deter2}) is thus equal to
\Be
\Delta = \zeta^{-12}\left(\frac{27\sqrt{3} a}{1024\sqrt{2}} \right)^3 \det\left(\begin{array}{ccccc}
\bft_2 -\bft_1 & \bft_3-\bft_1&\cdots &\bft_4-\bft_3\\
\bft_2 \cross\bft_1 & \bft_3\cross\bft_1&\cdots &\bft_4\cross\bft_3\
\end{array}\right)+O(\zeta^{-13})
\Ee
which does not vanish for $\zeta$ large enough.

\subsection{From local to global controllability}
\label{localtoglobal}
In all the three model swimmers that we proposed, we have found a point $(\xi_0,p_0)\in \mathcal{S}\times \mathcal{P}$ at which the sufficient condition for local controllability 
\begin{equation}
\mbox{ det}\left(\bfF_1,\cdots,\bfF_M,[\bfF_1,\bfF_2],\cdots,[\bfF_{M-1},\bfF_M]\right)(\xi_0,p_0)\ne 0\,,
\label{condcontrol}
\end{equation}
holds, which proves, denoting by $\mathcal{F}$ the family of analytic vector-fields $(\bfF_1,\cdots,\bfF_M)$, that
$$
\mbox{dim } Lie_{(\xi_0,p_0)}(\mathcal{F})=M+\frac{M(M-1)}{2}=d\,.
$$
In order to show the global controllability of the problem, we show the following lemma
\begin{lemma}
For every $(\xi,p)\in \mathcal{S}\times \mathcal{P}$, 
\Ben
\mbox{dim } Lie_{(\xi,p)}(\mathcal{F})=d\,.
\label{condglobal}
\Een
\label{lemmaglobal}
\end{lemma}

\begin{proof}
Let $(\xi,p)\in \mathcal{S}\times\mathcal{P}$. From Hermann-Nagano Theorem (see \cite{Jurdjevic-1997-a}, p. 48), we know that (\ref{condglobal}) holds at every point $(\bar{\xi},\bar{p})$ in the orbit $\mathcal{O}_{\mathcal{F}}(\xi_0,p_0)$ of $\mathcal{F}$ passing through $(\xi_0,p_0)$. Now, due to the specific form of the $(\bfF_i)_{1\leq i\leq M}$, there exists $q\in \mathcal{P}$ such that $(\xi,q)\in \mathcal{O}_{\mathcal{F}}(\xi_0,p_0)$. Indeed, 
$$
(\xi,q) = \exp\left(\sum_{i=1}^M (\xi-\xi_0)_i\bfF_i\right)(\xi_0,p_0)\,.
$$
Eventually, because of the invariance of the problem with respect to the orientation of the swimmer, one has
$$
\mbox{dim } Lie_{(\xi,p)}(\mathcal{F})=\mbox{dim } Lie_{(\xi,q)}(\mathcal{F})=\mbox{dim } Lie_{(\xi_0,p_0)}(\mathcal{F})=d\,.
$$
This ends the proof of Lemma \ref{lemmaglobal} and thus of Theorem \ref{ourtheorem}.
\end{proof}
\section{Optimal swimming}
\label{sec5}
A possible notion of optimal swimming consists in minimizing the
energy dissipated while trying to reach a given position $p_T$ at time
$T$ from an initial position $p_0$ at time 0. In this sense, the
optimal stroke is the one that minimizes
\begin{equation}
  \label{eq:energy-of-v}
  \mathcal{E}(\xi, p) := \int_0^T   \sum_{i=1}^N  \sum_{j=1}^N
  \int_{\pB} \bfv_i(\xi,p,\bfr)\cdot \TT^{-1}_{\bfx(\xi,p)}
 \bfv_j(\xi,p,\bfr)\,d\sigma(\bfr) \,dt
\end{equation}
where $\bfv$ satisfies the self-propulsion conditions expressed by
Equations~\eqref{noforce} and~\eqref{notorque}, subject to the
constraint that $p(T)$ is equal to $p_T$.

The special structure of our problem allows us to express the
dissipated energy solely in terms of the shape variables $\xi$, by
virtue of equation \eqref{linear-relationship-xi-dot-p}.  We start by
expressing the velocity of the $m-$th ball due to a unit shape change
in the $i$-th component of $\dot \xi$ as
\begin{equation}
  \label{eq:def-wi}
  \bfw^m_i(\xi, p, \bfr) :=  \nabla_{\xi} \bfX_m(\xi,p,\bfr) \cdot
  \bfe_i + \nabla_{p} \bfX_m(\xi,p,\bfr) \cdot (\VV(\xi, p) \cdot \bfe_i)\,.
\end{equation}

This allows us to simplify the energy of the system to an expression
depending on the rate of shape and position changes $\dot \xi$:
\begin{equation}
  \label{eq:energy-of-xi-p}
  \mathcal{E}(\xi, p) := \int_0^T   \sum_{i=1}^M  \sum_{j=1}^M
  \dot \xi_i \MetG_{ij}(\xi, p) \dot \xi_j\,dt,
\end{equation}
where the metric $\MetG$ is defined by
\begin{equation}
  \label{eq:metric-G}
  \MetG_{ij}(\xi, p) :=  \sum_{m=1}^N  \int_{\pB}
  \bfw_i^m(\xi, p, \bfr)
  \cdot \TT^{-1}_{\bfx(\xi,p)} (\bfw_j^1(\xi, p, \bfr),\cdots,\bfw_j^N(\xi, p, \bfr))^m \,d\sigma(\bfr).
\end{equation}

Since dissipation is invariant with respect to the current position
$p$, it is possible to show that $\MetG_{ij}(\xi,p) = \MetG_{ij}(\xi, 0) :=
\MetG_{ij}(\xi)$ for any position $p$.

Optimal swimming is then given by the solution $(\overline\xi, \overline p):[0,T]\longrightarrow \mathcal{S}\times\mathcal{P}$ of
\begin{subequations}
  \label{eq:optimization}
  \begin{align}
    \min_{(\xi,p) \in  (\mathcal{S}\times\mathcal{P})}\qquad & \int_0^T \dot \xi \cdot
    \MetG(\xi) \cdot \dot
    \xi \,d t =: \mathcal E(\xi) \label{eq:opt-minimality}\\
    \text{subject to }\qquad & C(\xi, p, t) = 0 && \forall t \in [0,T]\label{eq:opt-feasibility} \\    
    & \xi_L  \leq \xi(t) \leq \xi_U && \forall t \in [0,T],\label{eq:shape-bounds}
  \end{align}
\end{subequations}
where the contraints $C(\xi, p, t) = 0$ are given explicitly by
\begin{equation}
  \label{eq:constraints}
  \begin{aligned}
    & \VV(\xi(t),p(t))\cdot\dot \xi(t) - \dot p(t)
    & =  0,  && \quad \forall t \in [0,T],\\
    & \xi(0)  - \xi(T) & = 0,\\
    & p(0) - p_0  & = 0,\\
    & p(T) - p_T  & = 0,\\
    & \xi(0) -\xi_0 & = 0.
 \end{aligned}
\end{equation}

Some of the constraints in~\eqref{eq:constraints} can be relaxed. For
example, it is possible to fix only the periodicity of the shape, but
not its actual initial and final values (so that the last constraint
in~\eqref{eq:constraints} disappears), or to fix only some of the
components of $p(T)$, leaving the others free to be optimized.

\section{Numerical Approximation}
\label{sec6}
Our strategy for the solution to problem~\eqref{eq:optimization} is to
discretize in time both $\xi(t)$ and $p(t)$ using cubic splines,
and rewrite everything as a finite dimensional constrained
optimization problem, defined only on the coefficients of the spline
approximation of $\xi$ and $p$. 

A \texttt{C++} code for boundary integral methods, which has been developed using  the
\texttt{deal.II} library (see~\cite{BangerthHartmannKanschat-2007-a}
and~\cite{BangerthHartmannKanschat--a}), solves the Dirichlet-to-Neumann map $\mathcal{T}^{-1}$, and is then used to construct
$\MetG(\xi)$ and $\VV(\xi,p)$.

The finite element method in the spline space is then used to produce
a finite dimensional approximation to~\eqref{eq:optimization}, which
is solved using the reduced space successive quadratic programming
algorithm (rSQP) provided with the \texttt{Moocho} package of the
\texttt{Trilinos} \texttt{C++}
library~\cite{HerouxBartlettHowle-2005-a}.

In this section we touch only briefly on the numerical implementation
of the discrete Dirichlet-to-Neumann map for a collection of spheres,
and we refer to~\cite{AlougesDeSimoneHeltai-2011-a} for more
details on the rSQP technique used for the search of constrained
minima.

In principle, standard surface meshing techniques could be employed to
integrate on the surface of our swimmer, but they lack the speed and
accuracy that specialized methods can offer when integrating on such
special domains like collections of spheres.

The method we have chosen for the integration on the spheres is presented
in~\cite{HannayNye-2004-a} for use with smooth functions and
experimental data. It uses an oblique array of integration sampling
points based on a chosen pair of successive Fibonacci numbers.

The pattern which is obtained on each sphere has a familiar appearance
of intersecting spirals, avoiding the local anisotropy of a
conventional latitude-longitude array. We modified the original method
to take explicitly into account multiple spheres and the singularity
of the Stokeslet (see~\cite{Pozrikidis-1992-a} for a more detailed
introduction on boundary integral methods for viscous flow).

The discrete versions $\bfu_h$ and $\bff_h$ of the continuous velocity
and force densities $\bfu$ and $\bff$ are obtained by sampling the
original functions at the quadrature points, translated and replicated
on each sphere according to $\xi$ and $p$, to obtain $N_V$ scalars
that represent $N_V/3$ vectors applied to the quadrature points
sampled on $\pO$.

Integration on $\pO$ of the function $\bfu_h(\bfx)$ is then
effectively approximated by
\begin{equation*}
  \int_\pO \bfu(\bfx) \,d \sigma(\bfx) \sim \Francois{\sum_{i=1}^{N_V/3}}\bfu(\bfq^i) \omega^i,
\end{equation*}
where the weights $\omega^i$ and quadrature points $\bfq^i$ are
derived according to~\cite{HannayNye-2004-a}. 

The construction of a discrete Dirichlet-to-Neumann map for our
swimmers follows the usual collocation method
(see, for example,~\cite{Pozrikidis-2002-a}),
\begin{equation}
\label{eq:collocation-DN-map}
\begin{aligned}
  \bfu(\bfq^i) & = \int_\pO \bfG(\bfy-\bfq^i) \bff(\bfy)
  \,d\sigma(\bfy) &&\\
  & \sim \sum_{j=1, j\neq i}^{N_V/3}\bfG(\bfq^j-\bfq^i) \bff(\bfq^j) \omega^j +
  \overline \bfG^{ii} \bff(\bfq^i) \omega^i  &&\qquad i=\Francois{1}, \dots, N_V/3.
\end{aligned}
\end{equation}

A specialized calculation of the singular integral is performed in
Equation~\eqref{eq:collocation-DN-map} when $i$ is equal to $j$, by
deriving $\overline \bfG^{ii}$ according to the (known) exact solution
of the flow around a sphere.

Equation~\eqref{eq:collocation-DN-map} can be expressed more compactly
as:
\begin{equation}
  \label{eq:discrete-DN-map}
  \bff_h  = \mathcal A^{-1} \bfu_h,
\end{equation}
where the matrix $\mathcal A$ is obtained by  writing explicitly the
various components of Equation~\eqref{eq:collocation-DN-map}, and we
have used the notation $\bfu_h$ and $\bff_h$ to indicate the $N_V$
dimensional vector that samples the functions $\bfu(\bfx)$ and
$\bff(\bfx)$ at the quadrature points $\bfq$.

The construction of the discrete Dirichlet-to-Neumann map allows us to
compute explicitly, given a pair of shape and position variables $(\xi,p)$,
the numerical approximations of $G_h(\xi)$ and $V_h(\xi,p)$. These can in turn be
used to solve the constrained minimization problem as in~\cite{AlougesDeSimoneHeltai-2011-a}, i.e., by
using the implementation of the rSQP method included in the
\texttt{Moocho} package of the \texttt{Trilinos}
library~\cite{Bartlett-2009-a}.

\section{Numerical Results}
\label{sec:numerical-results}

In this section we analyze some of the optimal strokes that can be
obtained with our software for the case of the three-sphere swimmer in
the plane~\cite{Lefebvre-LepotMerlet-2009-a} (3SP), depicted in
Figure~\ref{3sub}. While this is one of simplest possible
extension of the three-sphere swimmer by Najafi and
Golestanian~\cite{NajafiGolestanian-2004-a}, it adds already a high
degree of complexity to the numerical approximation scheme.

The added complexity comes both from the presence of three additional time
dependent variables (one shape and two positions), and from the fact
that the vector field that relates  positional changes to shape
changes ($\dot p = \VV(\xi, p)\cdot \dot \xi$) depends also (though in an
explicit way) on the orientation $\theta = p_3$ of the swimmer.


In the software used to compute the optimal strokes for axisymmetric
swimmers~\cite{AlougesDeSimoneHeltai-2011-a}, there is only one
position variable $p_1(t) = c(t)$, which can be obtained by
post-processing the vector field $\VV(\xi)$ (which, in the axisymmetric
case, is independent on the position $p$) and $\xi$ itself.

In the present case of a planar swimmer, this is no longer possible, and the optimization
software has to go from optimizing over the \emph{two} shape variables
$\xi_1(t)$ and $\xi_2(t)$ to optimizing over the entire shape-position
space, composed of the \emph{six} time dependent variables
$\xi(t) = (\xi_1(t),\xi_2(t), \xi_3(t))$ and $p(t) = (c_x(t), c_y(t),
\theta(t))$.

In the experiments that follow, we use water at room temperature
($20^{o}$ C) as the surrounding fluid, and we express lengths in
millimeters ($mm$), time in seconds ($s$) and weight in milligrams
($mg$). Using these units, the viscosity of water is approximately one
($1 mPa\, s = 1 mg\, mm^{-1}\, s^{-1}$), and the energy is expressed in
pico-Joules ($1pJ=1mg\, mm^2\, s^{-2}=10^{-12} J$).

The radius of each ball of the swimmer is $.05mm$, and additional
constraints are added to the optimization software to ensure that each
component of the shape variable $\xi$ stays in the interval $[.1mm,
.7mm]$, so that the balls never touch each other and that they don't
separate too much.

We present some of the computed optimal strokes that illustrate the
richness of the problem. In all the examples that follow, we ask the
swimmer to go from $p(0) = (0mm, 0mm, \theta_0)$ to $p(T) = (.01mm,
0mm, \theta_T)$, for various combinations of $\theta_0$, $\theta_T$
and $\xi(0)=\xi(T)$.

Since the swimmer is $2\pi/3$ periodic, it is sufficient to study the
optimal strokes varying $\theta_0$ between $0$ and $2\pi/3$ while
keeping $c_x$ and $c_y$ constant to explore all possible optimal
swimming strokes. All other combinations can be obtained by
rotations and permutations between the components of the optimal
solution $(\xi, p)$.

Moreover, since rotating the swimmer by $\pi$ is equivalent to asking
it to swim in the opposite direction, reversibility of Stokes flow
allows us to further reduce the study of the optimal strokes to the
interval $\theta_0 \in [0,\pi/3]$.

In particular, if the pair $\xi(t)$ and $p(t) = (\bfc(t) , \theta(t))$
is an optimal solution of \eqref{eq:optimization}, then we have that
all the following permutations, symmetries and rotations of $(\xi,p)$
are also optimal swimming strokes, but with different initial and
final positions:
\begin{subequations}
\label{eq:permutations}
\begin{alignat}{2}
  \xi(t)\qquad & \bfc(t)+k, \theta(t) &&
  \forall k \in \Re^2 \label{eq:perm-add}\\
  \xi(t)\qquad & \Rot_z(\gamma)\bfc(t),
  \theta(t)
  + \gamma && \forall \gamma \in [0, 2\pi] \label{eq:perm-rot}\\
  \xi_2(t), \xi_3(t), \xi_1(t) \qquad &
  \Rot_z(2\pi/3)\bfc(t),\theta(t) && \label{eq:perm-2/3pi}\\
  \xi_3(t), \xi_1(t), \xi_2(t) \qquad &
  \Rot_z(4\pi/3)\bfc(t),\theta(t) && \label{eq:perm-4/3pi}\\
  \xi(T-t) \qquad & \bfc(T-t)-\bfc(T), \theta(T-t)\label{eq:perm-inv}.
  &&
\end{alignat}
\end{subequations}

Given a family of optimal solutions with different starting angles
$\theta_0 \in [0,\pi/3]$, then properties~\eqref{eq:permutations} can be
used to derive the optimal solutions in all of the interval $[0,2\pi]$.

The optimal strokes we present have been formatted to show the
trajectory of the center of the swimmer amplified by a factor 60 with
respect to the dimensions of the swimmer itself, and the orientation
$\theta$ of the swimmer amplified by a factor 10. The full
configuration of the swimmer is shown at ten equally spaced times between $t=0$ and $t=T=1$.

\begin{figure}[!htb]
  \centering
  
  \resizebox{.5\textwidth}{!}{
    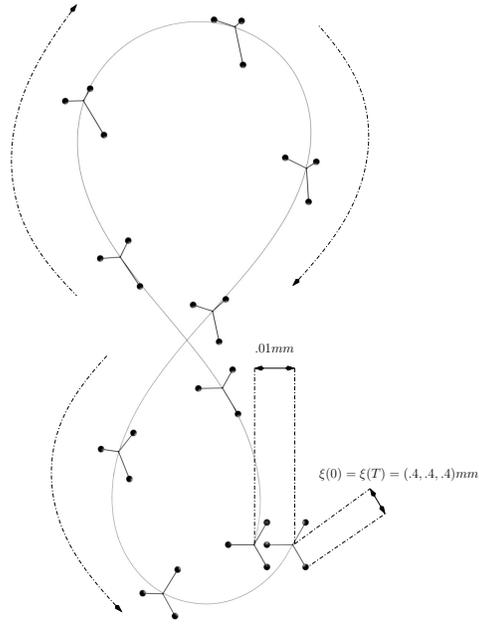
  }

  \caption{3SP swimmer: optimal stroke for $\theta_0=0=\theta_T$
    radians. The initial shape is fixed at $\xi(0) = \xi(T) = (.4, .4,
    .4)mm$. Energy consumption to swim by $.01mm$ in the $x-$direction:
    $0.511pJ$.}
  \label{fig:plane0}
\end{figure}

Our first optimal stroke (Figure~\ref{fig:plane0}) is obtained by
fixing the initial and final rotation $\theta_0 = \theta_T$ to be zero
and the initial shape of the swimmer to be $\xi = (.4, .4, .4)mm$. The
total energy consumption to swim by $.01mm$ in the $x$ direction with a
single stroke is $0.511pJ$.

The first important thing to notice is that, although the initial
configuration is symmetric with respect to the translation axis, the
stroke is not. This implies that the optimal solution is not unique:
exchanging the off axis shape variables $\xi_2(t)$ and $\xi_3(t)$
yields to the same final displacement with the same total energy
consumption. 

A configuration which is symmetric with respect the $y-$axis can be
obtained by setting $\theta_0 = \pi/6$, and also in this case the
solution is non-symmetric and non-unique.

Given the permutation properties~\eqref{eq:permutations}, we deduce
that two distinct global solutions with the same energy dissipation
exist for $\theta_0 = N\pi/6$, for any integer $N$. From the existence
of these different global minima we can deduce that at least two
branches of \emph{local} minima exist outside those particular points,
starting off as continuous perturbations of the two global solutions.

\begin{figure}[!tb]
  \centering
  
  \resizebox{.5\textwidth}{!}{
    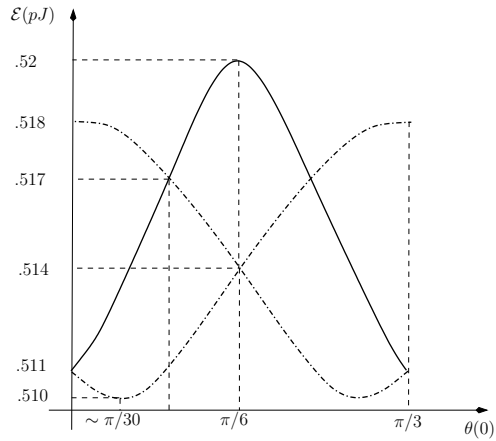
  }

  \caption{Branches of optimal swimming for 3SP swimmer: dissipated energy VS starting
    angle.}
  \label{fig:branches}
\end{figure}

Numerical evidence show that at least three stable branches exist, as
shown in Figure~\ref{fig:branches}, where we plot the energy
dissipation of all the (local) optimal strokes with respect to the
starting initial angle $p_3(0) = \theta_0$. It is possible to steer
the software towards one or the other local branch by feeding it with
initial guesses which are close to the desired branch.

From Figure~\ref{fig:branches} it is also clear that the swimmer
prefers to swim along directions which are close to parallel to one of
the arms: $N\pi/3\pm\pi/30$ are the approximate preferred swimming
directions. The directions perpendicular to one of the arms ($\theta_0
= \pi/6+N\pi/3$) are the most expensive ones (although they differ in energy only by approximately 2\%).

\begin{figure}[!htb]
  \centering
  
  \resizebox{.7\textwidth}{!}{
    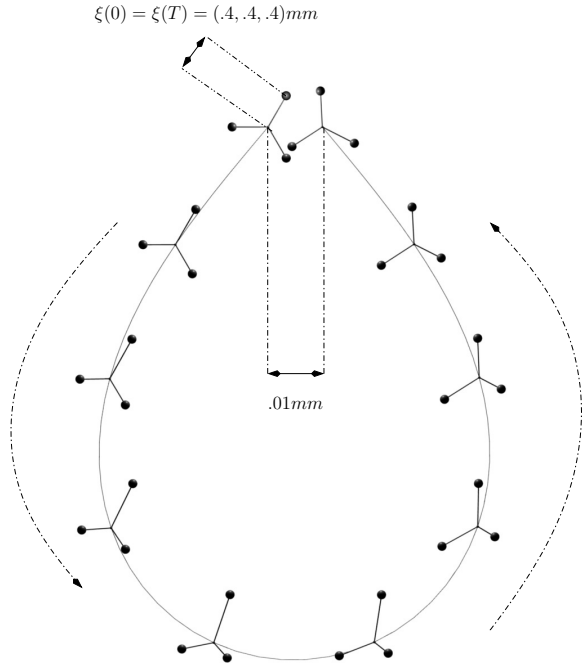
  }

  \caption{3SP swimmer: optimal stroke for $\theta_0=0$, $\theta_T$
    optimized to $.058$ radians. The initial shape is fixed at $\xi(0)
    = \xi(T) = (.4, .4, .4)mm$. Energy consumption to swim by $.01mm$ in
    the $x-$direction: $.209pJ$.}
  \label{fig:plane1}
\end{figure}

From Figure~\ref{fig:plane0} it is evident that a lot of energy is
wasted in forcing the final rotation to be equal to the initial
one. If we relax this requirement and allow the software to optimize
the stroke without forcing the rotation to be periodic in time, then
the swimmer does not follow an ``eight-shaped'' path anymore, and the
actual target displacement is reached with a smaller stroke
(Figure~\ref{fig:plane1}), and by dissipating less than half the
energy ($.209pJ$ instead of $.511pJ$). The final rotation is
approximately $.058$ radians (about $3.3$ degrees). 

The main defect associated with this strategy is that all successive
strokes will be rotated with respect to the initial swimming
direction. If the stroke were to be repeated, we would observe a
rotational drift and an average circular motion along a trajectory of
radius $r=.35mm$. To ensure that \Francois{several} successive strokes
produce an average trajectory along a straight line, one cannot use
copies of the same stroke.  Instead, different optimization problems
need to be solved \Francois{for each stroke}, in each of which the
initial orientation is the one reached at the end of the previous
stroke.

\begin{figure}[!htb]
  \centering
  
  \resizebox{.8\textwidth}{!}{
    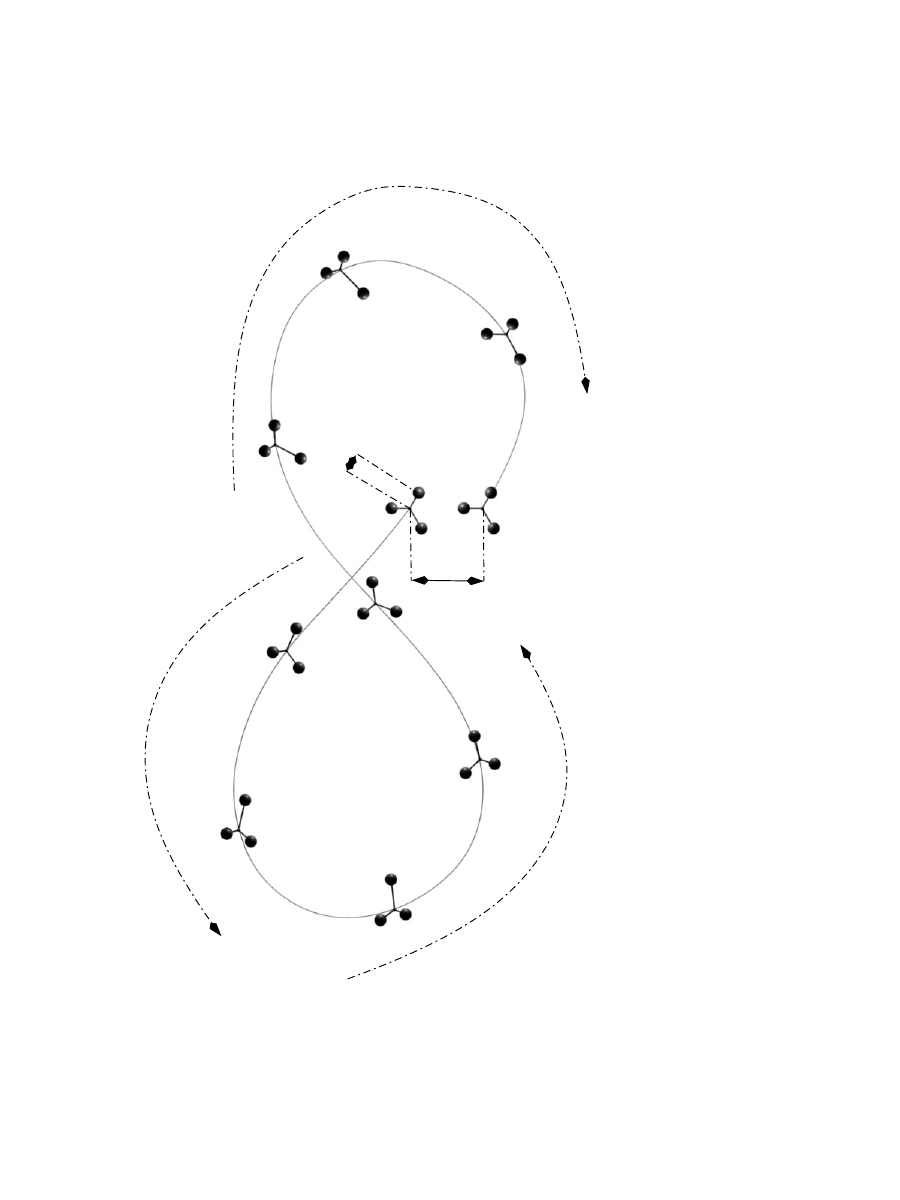
  }

  \caption{3SP swimmer: optimal stroke for $\theta_0=0 = \theta_T$.
    The initial shape is optimized to $\xi(0) = \xi(T) = (.15, .15,
    .19)mm$. Energy consumption to swim by $.01mm$ in the $x-$direction:
    $.221pJ$.}
  \label{fig:plane3}
\end{figure}

A complementary strategy to reduce the energy of a single stroke
without generating a rotated final configuration consists in
optimizing also on the starting shape, rather than on the final
angle. Figure~\ref{fig:plane3} shows this strategy, where the final
angle is fixed to be equal to the initial angle. From the figure it is
clear that the more complex hydrodynamic interaction due to the
closeness of the spheres can be effectively exploited by the
optimization software to reduce the dissipated energy.

The optimal initial shape is found to be $(.15, .15, .19)mm$, a non-symmetric shape. The total energy dissipation with this stroke
strategy is about $.221pJ$. While this is already less than half of
the energy consumed by our first optimal stroke
(Figure~\ref{fig:plane0}), it is still about 10\% more expensive than
the second attempt (Figure~\ref{fig:plane1}).

\begin{figure}[!htb]
  \centering
  
  \resizebox{.6\textwidth}{!}{
    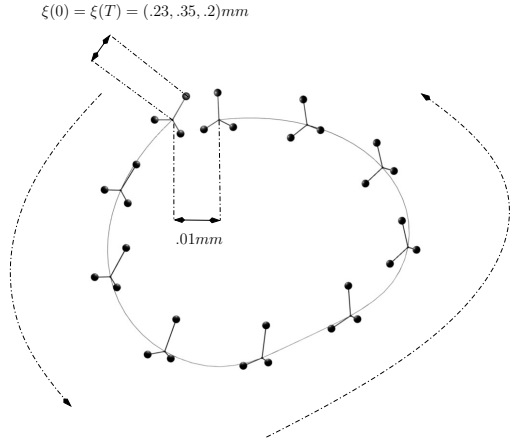
  }

  \caption{3SP swimmer: optimal stroke for $\theta_0=0$, $\theta_T$
    optimized to $.057$ radians.  The initial shape is optimized to
    $\xi(0) = \xi(T) = (.23, .35, .2)mm$. Energy consumption to swim by
    $.01mm$ in the $x-$direction: $.100pJ$.}
  \label{fig:plane2}
\end{figure}

The combination of the two strategies, in which both the initial shape
and the final angle are left free for the software to optimize, is
shown in Figure~\ref{fig:plane2}. In this case the total energy
dissipation drops down to $.100pJ$. The initial shape is optimized to
the value $(.23, .35, .2)mm$ and the final angle is $0.057$ radians, or about
$3.28$ degrees. In this case the optimal stroke hits the boundary of
the region of allowed shapes.

{\bf{Acknowledgments:}} The first  and last authors benefited from the support of the 
``Chair Mathematical modelling and numerical simulation,
  F-EADS -- Ecole Polytechnique -- INRIA -- F-X''

\bibliographystyle{plain}
\bibliography{parking}

\end{document}